\documentclass[10pt]{amsart}
\usepackage{amsmath,amsfonts,amssymb,amsthm,amscd}
\newtheorem{rem}{Remark}
\newtheorem{thm}{Theorem}
\newtheorem{lem}[thm]{Lemma}
\newtheorem{prop}[thm]{Proposition}
\newtheorem{defn}[thm]{Definition}

\begin{document}

\title{The Bochner--flat cone of a CR manifold}
\date{}
\author{Liana David}
\maketitle

\textbf{Abstract:} We construct a K\"{a}hler structure (which we
call a generalised K\"{a}hler cone) on an open subset of the cone
of a strongly pseudo-convex CR manifold endowed with a
$1$-parameter family of compatible Sasaki structures. We determine
those generalised K\"{a}hler cones which are Bochner-flat and we
study their local geometry. We prove that any Bochner-flat
K\"{a}hler manifold of complex dimension bigger than two is
locally isomorphic to a generalised K\"{a}hler cone.

\section{Introduction}

The Bochner tensor of a K\"{a}hler manifold is the biggest
irreducible component of the curvature tensor under the action of
the unitary group. In complex dimension two, the Bochner tensor
coincides with the anti-self dual Weyl tensor. A K\"{a}hler
manifold is Bochner-flat if its Bochner tensor vanishes.
Bochner-flat K\"{a}hler manifolds represent an important class of
K\"{a}hler manifolds and have been intensively studied: the local
geometry of Bochner-flat K\"{a}hler manifolds and its interations
with Sasaki geometry has been studied, using the Webster's
correspondence, in \cite{paul}; complete Bochner-flat K\"{a}hler
structures on simply connected manifolds have been classified in
\cite{bryant}; generalisations of Bochner-flat K\"{a}hler
manifolds (like weakly Bochner-flat K\"{a}hler manifolds and
K\"{a}hler manifolds with a hamiltonian $2$-form) have also been
developed (see, for example, \cite{acg}, \cite{acg2}, \cite{tim}).

An important class of K\"{a}hler manifolds is represented by the
K\"{a}hler cones of Sasaki manifolds. Unfortunately, except when
the Sasaki manifold is an open subset of the standard CR sphere
with its standard metric as the Sasaki metric, the K\"{a}hler
cones are not Bochner-flat. In this paper we propose an
alternative construction, which is a natural generalisation of the
K\"{a}hler cone construction and which produces, locally, all
Bochner-flat K\"{a}hler structures of complex dimension bigger
than two. More precisely, we consider, on a fixed CR manifold
$(N,H,I)$, a $1$-parameter family of Sasaki Reeb vector fields $\{
T_{r}, r\in \mathcal J\}$ (with $\mathcal J\subset\mathbb{R}^{>0}$
an open connected interval). On the cone manifold $N\times\mathcal
J$ we define an almost complex structure $J$, which on $H\subset
T(N\times{\mathcal J})$ coincides with $I$ and which sends the
radial vector field $V=r\frac{\partial}{\partial r}$ to the vector
field $T$, defined by $T(p,r):=T_{r}(p)$, for any $(p,r)\in
N\times\mathcal J .$  It turns out that $J$ is integrable and that
the pair $(\omega :=\frac{1}{4}dd^{J}r^{2},J)$ is a K\"{a}hler
structure on the open subset of $N\times\mathcal J$ where $\omega
(V,T)>0.$ Such a K\"{a}hler structure will be  called a
generalised K\"{a}hler cone and coincides with the usual
K\"{a}hler cone of a Sasaki manifold when the family of Reeb
vector fields is constant. A strong motivation for this
construction comes from the fact that the Bryant's family of
Bochner-flat K\"{a}hler structures (which have been discovered by
Robert Bryant in his classification theorem of complete
Bochner-flat K\"{a}hler structures on simply connected manifolds
\cite{bryant} and have been further studied in \cite{paul}) are
generalised K\"{a}hler cones. Our main result is the following:

\begin{thm}\label{main}
Any Bochner-flat K\"{a}hler manifold of complex dimension bigger
than two is locally isomorphic to a generalised K\"{a}hler cone.\end{thm}

The plan of the paper is the following: in Section \ref{basic} we
review the theory of K\"{a}hler and Sasaki manifolds, which will
be useful later on in our study of generalised K\"{a}hler cones.
In Section \ref{det} we determine the generalised K\"{a}hler cones
which are Bochner-flat and in Section \ref{quotient} we study
their local geometry. This study will readily imply Theorem
\ref{main}. The last Section is dedicated to examples. We explain
how K\"{a}hler manifolds with constant holomorphic sectional
curvature, weighted projective spaces and Bryant's family of
Bochner-flat K\"{a}hler structures fit into our formalism of
generalised K\"{a}hler cones. We also look at Bochner-flat
K\"{a}hler generalised K\"{a}hler cones of order one and at those
which are of Tachibana and Liu type
\cite{liu}.\\

\textbf{Acknowledgements:} I warmly thank Paul Gauduchon and David
Calderbank for encouraging me to study Bochner-flat generalised
K\"{a}hler cones and for many useful comments on a first version
of this paper. I thank the referee for 
his detailed report and his
deep observations and suggestions.

\textit{Key words and phrases:}
Bochner-flat K\"{a}hler manifolds, cone manifolds, CR and Sasaki manifolds.

\textit{Mathematics Subject Classification:}
53C55, 53C21, 53D10.

\section{Notations and earlier results}\label{basic}

\subsection{The Bochner tensor of a K\"{a}hler manifold}
In this section we recall the definition of the Bochner tensor
of a K\"{a}hler manifold. We use the formalism developed in
\cite{acg}, \cite{tim}.

Let $(V,g,J)$ be a real vector space together with a complex
structure $J$ and a $J$-invariant positive definite
metric $g$.
We shall identify vectors and covectors of $V$ using the metric $g$.
Let $\omega :=g(J\cdot ,\cdot )$ be the K\"{a}hler form.
Recall that the space
${\mathcal K}(V)$ of K\"{a}hler curvature
tensors of $(V,g,J)$, defined as those
curvature tensors which annihilate all $J$-anti-invariant $2$-forms
on $V$, decomposes into a $g$-orthogonal sum
\begin{equation}\label{decompo}
{\mathcal K}(V):=c_{\mathcal K}^{*}\left(\mathrm{Sym}^{1,1}(V)\right)\oplus
{\mathcal W}(V),
\end{equation}
where $c_{\mathcal K}^{*}:\mathrm{Sym}^{1,1}(V)\to{\mathcal K}(V)$
is the adjoint of the Ricci contraction
$$
c_{\mathcal K}:{\mathcal K}(V)\to \mathrm{Sym}^{1,1}(V), \quad
c_{\mathcal K}(R)(v,w):=\mathrm{trace}R(v,\cdot ,w,\cdot ), \quad v,w\in V
$$
and has the following expression \cite{tim}
\begin{equation}\label{dual}
c^{*}_{\mathcal K}(S)=\frac{1}{2}[\frac{S\land
\mathrm{Id}+(J\circ S)\land J}{2}+\omega\otimes S
+\beta\otimes J],
\end{equation}
where $S\in\mathrm{Sym}^{1,1}(V)$ is a symmetric $J$-invariant
endomorphism of $V$, "$\mathrm{Id}$" is the identity
endomorhism,
$\beta\in
\Lambda^{1,1}(V)$ is the $J$-invariant $2$-form on $V$,
related to $S$ by
$\beta (v,w):=g(SJv,w)$, and, for two endomorphisms $S$ and $T$
of $V$, $S\land T$ is the endomorphism of $\Lambda^{2}(V)$ defined
by the formula
$$
(S\land T)(v,w):=S(v)\land T(w)-S(w)\land T(v),\quad v,w\in V.
$$
According to the decomposition (\ref{decompo}), a K\"{a}hler curvature tensor
$R\in{\mathcal K}(V)$ decomposes into the sum
$$
R=c_{\mathcal K}^{*}(S)+W^{K},
$$
where $W^{K}\in{\mathcal W}(V)$ is its principal part
(or the Bochner tensor of $R$)
and $S\in\mathrm{Sym}^{1,1}(V)$ is a modified Ricci tensor.

Consider now a K\"{a}hler manifold
$(M,g,J)$. The curvature $R^{g}$ of the K\"{a}hler
metric $g$, is, at every point $p\in M$, a K\"{a}hler curvature tensor
of the tangent space $(T_{p}M,g_{p},J_{p}).$
The principal
part of $R^{g}$ is called the Bochner tensor of $(M,g,J)$ and
is a section of
the symmetric product $\Lambda^{1,1}(M)\odot\Lambda^{1,1}(M)$.
The K\"{a}hler manifold
$(M,g,J)$ is Bochner-flat if its Bochner tensor vanishes.

\subsection{Review of CR and Sasaki manifolds}

Recall that an oriented (strongly pseudo-convex) CR manifold $(N,H,I)$ has a
codimension one oriented subbundle $H$ of the tangent bundle
$TN$, called the contact bundle, and a bundle homomorphism $I:H\to H$
with $I^{2}=-\mathrm{Id}$, such that, for every smooth sections
$X,Y\in\Gamma (H)$,
$[IX,IY]-[X,Y]$ is also a section of $H$ and
the integrability condition
\begin{equation}\label{inte}
[IX,IY]-[X,Y]=I\left( [IX,Y]+[X,IY]\right)
\end{equation}
is satisfied. Since $N$ and $H$ are
oriented, the co-contact line bundle $L:=TN/H$ is also oriented,
hence trivialisable.
A positive section $\mu$ of $L$ defines a contact form
$\theta :=\eta \mu^{-1}$ on $M$, where $\eta :TN\to L$ is
the natural projection and $\mu^{-1}\in \Gamma (L^{*})$
is the dual section of $\mu$, i.e.
the natural contraction between $\mu$ and $\mu^{-1}$ is the function
on $N$ identically equal to one. The bilinear form
$g(X,Y):=\omega(X,IY):=\frac{1}{2}d\theta (X,IY)$ of the bundle
$H$ is independent, up to a positive multiplicative function,
of the choice of the contact form and is positive definite - the strongly
pseudo-convexity condition. The contact form $\theta$ determines a
Reeb vector field $T$, uniquely
defined by the conditions $\theta (T)=1$
and $i_{T}d\theta =0$. Note that the Reeb vector field
preserves the bundle $H$, i.e. $[T,X]\in\Gamma (H)$ when $X\in\Gamma (H)$
and hence ${\mathcal L}_{T}(I)$ is a well-defined endomorphism of $H$.
There is also a Riemannian
metric $g$ of $N$ associated to $\theta$, which on $H$ is defined above
and such that $T$ is of norm one and
orthogonal to $H$. Finally, we need
to recall the definition of the Tanaka connection \cite{tanaka} associated
to $\theta$. It is the unique connection
$\nabla$ on $N$ with the following three properties: (i) it preserves
the bundle $H$; (ii) $I$, $g$ and
$T$ are $\nabla$ parallel; (iii) the torsion $T^{\nabla}$ of $\nabla$
has the following expression:
\begin{align*}
T^{\nabla}(X,Y)&=2\omega (X,Y)\\
T^{\nabla}(T,X)&=-\frac{1}{2}I{\mathcal L}_{T}(I)(X),
\end{align*}
for every  $X,Y\in\Gamma (H).$
It turns out that on $H$, $\nabla$ is determined
by a Koszul type formula
\begin{align*}
2g(\nabla_{X}Y,Z)&=X\left( g(Y,Z)\right) +Y\left( g(X,Z)\right) -
Z\left( g(X,Y)\right)\\
&+g([X,Y]^{H},Z)-g([X,Z]^{H},Y)-g([Y,Z]^{H},X)\\
\end{align*}
where $X,Y,Z\in\Gamma (H)$ and for a vector field $W$ of $N$,
$W^{H}:=W-\theta (W)T$ is its $g$-orthogonal projection on the bundle $H$.
The metric $g$ is called Sasaki
if, by definition, $T$ is a Killing vector field for the metric $g$, or,
equivalently, if
${\mathcal L}_{T}(I)=0$. In this case,
the curvature $R^{\nabla}$ of the Tanaka connection on the bundle $H$
is an element of the tensor product $\Lambda^{2}(N)\otimes\Lambda^{1,1}(H)$
and its restriction to bivectors of $H$
belongs to $\Lambda^{1,1}(H)\odot\Lambda^{1,1}(H)$
and is a K\"{a}hler curvature tensor of the
complex Riemannian vector bundle $(H,g,I)$. Its Bochner part --
called the
Chern-Moser tensor \cite{chern}, \cite{liana} of the CR manifold $(N,H,I)$ --
is independent of the choice of the compatible Sasaki
structure on $(N,H,I)$.
A CR manifold with vanishing
Chern-Moser tensor is called flat.
The importance of the Chern-Moser tensor comes from the fact that
if the CR manifold is flat and of dimension
bigger than three, then it is locally isomorphic with a sphere with
its standard CR structure \cite{chern}, \cite{vesti}.
On the other hand, if $g$ is Sasaki, the
complex structure $I$ and the metric $g$ of the bundle
$H$ descend
on the quotient $N/T$ and determine
a K\"{a}hler structure on this quotient.
(In our conventions, the quotient $N/T$ denotes the space of leaves
of the foliation generated by $T$ in a sufficiently small
open subset of $N$, so that $N/T$ is a manifold).
Moreover, the Bochner tensor of the K\"{a}hler manifold $N/T$ becomes
identified
with the Chern-Moser tensor of the CR manifold
$(N,H,I)$ \cite{web}, \cite{liana}.
Since any K\"{a}hler manifold
can be locally written as a quotient of a Sasaki manifold under its
Reeb vector field, by means of a choice of a local primitive
of the K\"{a}hler form, it follows that a Bochner-flat K\"{a}hler manifold
of complex dimension $m\geq2$, is locally isomorphic with
the quotient of a standard
CR sphere $S^{2m+1}$ under the Reeb vector field of a compatible
Sasaki structure.

\subsection{The local type of Bochner-flat K\"{a}hler structures}

The local geometry of Bochner-flat K\"{a}hler structures, of
complex dimension $m\geq 2$, is of four types \cite{bryant}, \cite{david},
\cite{paul}.
This follows from the fact that the compatible Sasaki structures
on the CR sphere
$S^{2m+1}$
are determined by elements of the Lie algebra $su(m+1,1)$,
and that there
are four conjugacy classes in this Lie algebra
(elliptic, hyperbolic, $1$-step and $2$-step  parabolic).
In order to explain this, it is convenient to identify $S^{2m+1}$
with the hermitian sphere $\Sigma^{2m+1}$
of all complex null
lines in a hermitian complex vector space $W$ of signature
$(m+1,1)$, with hermitian metric $(\cdot ,\cdot )$, by
fixing an orthonormal basis of $W$, i.e. a basis $\{ e_{0},\cdots ,
e_{m+1}\}$ with $(e_{0},e_{0})=-1$, $(e_{j},e_{j})=1$, for
$j\in\overline{1,m+1}$ and $(e_{i},e_{j})=0$ for $i\neq j$, and
associating
to a complex null line $x$ of $W$ its unique representative
of the form $e_{0}+u$, where $u\in S^{2m+1}$ belongs
to the unit sphere of the positive definite hermitian
vector space $\mathrm{Span}\{ e_{1},\cdots ,e_{m+1}\} .$
Let $\eta$ be the natural (line bundle valued)
contact form of $\Sigma^{2m+1}$:
$$
\eta (X):=\mathrm{Im}(\hat{X}w,w),
\quad X\in T_{x}\Sigma^{2m+1} ,\quad 0\neq w\in x,\quad
x\in\Sigma^{2m+1} ,
$$
where $\hat{X}\in\mathrm{Hom}_{\mathbb{C}}(x,W)$ is a representative
of $X\in \mathrm{Hom}_{\mathbb{C}}(x,W/x).$
A hermitian trace-free endomorphism $A$ of $W$ determines a Reeb
vector field $T^{A}$ of a Sasaki structure on the open subset
$$
\Sigma^{2m+1}_{A}:=\{ x\in\Sigma^{2m+1}:\quad (Aw,w)>0,
\quad w\in x,\quad w\neq 0 \} ,
$$
defined in the following
way: at a point $x\in\Sigma^{2m+1}_{A}$,
$T^{A}_{x}\in\mathrm{Hom}_{\mathbb{C}}(x,W/x)$ associates to a non-zero
vector $w\in x$, the
class of $iAw$ in $W/x.$ The contact form of $T^{A}$ is
$\eta_{A}:=\frac{\eta}{(Aw,w)}$, i.e.
$$
\eta_{A}(X)
=\frac{\mathrm{Im}(\hat{X}w,w)}{(Aw,w)},\quad X\in T_{x}\Sigma^{2m+1},\quad
0\neq w\in x,\quad x\in\Sigma^{2m+1}_{A}.
$$
Employing the notations of \cite{paul}, we shall denote by $M_{A}$
the induced K\"{a}hler structure on the quotient
$\Sigma^{2m+1}_{A}/T^{A}$ 
and by $g_{A}$ its K\"{a}hler metric.

We end this section with a simple Lemma on hermitian operators
which will play an important role in our treatment. For
completeness of the exposition, we include its proof.

\begin{lem}\label{identitate} Let
$A:W\to W$ be a
hermitian operator
on a complex vector space $W$ with a
hermitian
metric $(\cdot ,\cdot )$ of signature $(m+1,1)$. Suppose that
$A$ satisfies $(Aw,w)=0$, for any null
vector $w$ which belongs to a non-empty open subset $D$ of
$W.$ Then $A=\lambda\mathrm{Id}$, for
$\lambda\in\mathbb{R}$. If, moreover, $A$ is trace-free, then $A=0.$
\end{lem}

\begin{proof} Let $w=w_{t}$ be a curve in $D$,
with $w_{t}$ null for any $t$, $w_{0}=w\in D$
and $\dot{w}_{0}=X.$ Taking the derivative at $t=0$ of the equality
$(Aw_{t},w_{t})=0$ and using the fact that $A$ is hermitian, we get
${\mathrm{Re}}(Aw,X)=0$. In particular, we deduce that
$(Aw,X)=0$, for any null vector $w\in D$ and any $X\in W$, which is
hermitian orthogonal to $w$.
This implies that $Aw=\lambda w$, where $\lambda\in\mathbb{R}$ depends a priori on $w$.
It
follows that the map
\begin{equation}\label{wedge}
W\ni u\to Au\land
u\in\Lambda^{2}(W)
\end{equation}
vanishes when $u\in D$ is null. Being holomorphic, the map
(\ref{wedge}) must be identically zero. We deduce that for any
$u\in W$, $Au$ is a multiple of $u$ which clearly implies the
first claim. The second claim is trivial.
\end{proof}

\subsection{The Bryant minimal and characteristic polynomials}\label{brymin}

The local type of a Bochner-flat K\"{a}hler manifold $(M,g,J)$ is
encoded into the Bryant's minimal and characteristic polynomials,
which can be defined as follows. Let $S$ be the modified Ricci
operator which satisfies $c^{*}_{\mathcal K}(S)=R^{g}$ 
(where $R^{g}$ is the curvature of $g$)
and $P(t)$
be the characteristic polynomial of a new modified Ricci operator $\Theta$,
related to $S$ by
\begin{equation}\label{theta}
\Theta :=\frac{1}{4}\left( S-
\frac{\mathrm{trace}_{\mathbb{R}}(S)}{2(m+2)}\mathrm{Id}\right) ,
\end{equation}
where $m$ is the complex dimension of $M$. The Ricci operator
$\Theta$ has been introduced by Robert Bryant in \cite{bryant}. It
will be considered as a complex linear operator on the complex
vector bundle $(TM,J).$ Its trace is called the modified scalar
curvature of $(M,g,J).$ Denote by $\xi_{1},\cdots ,\xi_{l}$ the
non-constant roots of $P$ and by $P_{n}$ its non-constant part,
defined by $P_{n}(t):=(t-\xi_{1})\cdots (t-\xi_{l})$. The number
$l$ is called the order of $(M,g,J)$. On a dense open subset $M^{0}$ of
$M$, the eigenvalues $\xi_{j}$ (for any $j\in \{ 1,\cdots ,l\}$)
are simple, different from each other at any point and different,
at any point, from any constant eigenvalue of $\Theta$; the
functions $\xi_{1},\cdots ,\xi_{l}$ are functionally independent on $M^{0}$
and
\begin{equation}\label{grad}
|\mathrm{grad}_{g}(\xi_{j})|^{2} =
-4\frac{p_{m}(\xi_{j})}{P_{n}^{\prime}(\xi_{j})},\quad j\in\{
1,\cdots ,l\}
\end{equation}
where $p_{m}$ is a monic polynomial of degree $l+2$, with constant
coefficients, independent of $j$, called the Bryant minimal
polynomial of $(M,g,J)$. The Bryant characteristic polynomial
$p_{c}$ of $(M,g,J)$ is by definition the product of $p_{m}$ with
the constant part $P_{c}:=P/P_{n}$ of $P.$\\

Suppose now that $(M,g,J)\cong M_{A}=\Sigma^{2m+1}_{A}/T^{A}$, for
a hermitian operator $A$ of $W$. Denote by $\tilde{a}$ the reduced
adjoint operator of $A$, defined by
$$
\tilde{a}(t)=t^{l+2}+a_{1}t^{l+1}+\cdots +a_{l+2},
$$
where
$$
a_{k}:=A^{k}-\sigma_{1}(q_{A})A^{k-1}+\cdots
+(-1)^{k}\sigma_{k}(q_{A}),
$$
and $\sigma_{k}(q_{A})$ is the $k$ elementary symmetric function
of the roots of the minimal polynomial $q_{A}$ of $A$. The reduced
adjoint operator $\tilde{a}$ satisfies
$(tI-A)\tilde{a}(t)=q_{A}(t)I$, for any $t\in\mathbb{R}.$ The
eigenspace of $\Theta$ corresponding to a non-constant eigenvalue
$\xi_{j}$ is spanned by the vector field $L_{j}$ which, viewed as
a section of $H$, is equal, at a point $x\in\Sigma^{2m+1}_{A}$, to
$$
L_{j}(w)=\tilde{a}(\xi_{j})w\quad\mathrm{mod}w.
$$
The non-constant part $P_{n}(t)$ of the modified Ricci operator
$\Theta$ of $(M,g,J)$, viewed as a polynomial with function
coefficients defined on $\Sigma^{2m+1}_{A}$, is equal, at a point
$x\in \Sigma^{2m+1}_{A}$, to
$$
p_{A,x}(t):=\frac{(\tilde{a}(t)w,w)}{(Aw,w)},\quad w\in x,\quad
w\neq 0.
$$
The constants eigenvalues of $\Theta$ can also be described in
terms of $A$: if $\lambda$ is a constant eigenvalue of $\Theta$,
of multiplicity $n$, then it is a multiple eigenvalue of $A$, of
multiplicity $n+1$, and the eigenspace of $\Theta$, at a point
$x\in\Sigma^{2m+1}_{A}$, corresponding to $\lambda$ can be
identified with the intersection of the hermitian orthogonal
$x^{\perp}\subset W$ with the eigenspace of $A$ corresponding to
$\lambda$ (see \cite{paul}). The Bryant minimal and characteristic polynomials
$p_{m}$ and $p_{c}$ coincide with the minimal polynomial $q_{A}$,
respectively to the characteristic polynomial $Q_{A}$ of $A$
\cite{david}, \cite{paul}. The modified scalar curvature of
$M_{A}$, viewed as a function on $\Sigma^{2m+1}_{A}$, is equal, at
$x\in\Sigma^{2m+1}_{A}$, to
$-\frac{(A^{2}w,w)}{(Aw,w)}$, where $w\in x$ is non-zero.\\

The following Lemma will be useful in our study of generalised
K\"{a}hler cones and is an easy consequence of the theory
developed in \cite{paul}. For completeness of the exposition, we
include its proof.

\begin{lem}\label{aditional}
For any $t\in\mathbb{R}$, $x\in\Sigma^{2m+1}_{A}$ and $w\in x$ non-zero,
$$
g_{A}(dp_{A,x}(t),dp_{A,x}(t))=4\left(
q^{\prime}_{A}(t)p_{A,x}(t)-q_{A}(t)p_{A,x}^{\prime}(t)
-2tp_{A,x}^{2}(t) +p_{A,x}^{2}(t)\frac{(A^{2}w,w)}{(Aw,w)}\right).
$$
\end{lem}

\begin{proof} Via the metric
$g_{A}$, the $1$-form $dp_{A}(t)$ corresponds to the vector field
$L_{t}$, which, viewed as a section of $H$, is equal, at
$x\in\Sigma^{2m+1}_{A}$, to the endomorphism 
$$
L_{t}(w):=2\left( \tilde{a}(t)w -p_{A,x}(t)Aw\right)\quad
\mathrm{mod}w.
$$
Its square norm is equal to
\begin{align*}
g_{A}(L_{t},L_{t}) &=4\frac{(\tilde{a}(t)w,\tilde{a}(t)w)-
2p_{A,x}(t)(\tilde{a}(t)w,Aw)
+p_{A,x}^{2}(t)(A^{2}w,w)}{(Aw,w)}\\
&=4\left( q_{A}^{\prime}(t)p_{A,x}(t)-q_{A}(t)p_{A,x}^{\prime}(t)
-2tp_{A,x}^{2}(t) +\frac{p_{A,x}^{2}(t)(A^{2}w,w)}{(Aw,w)}\right),
\end{align*}
where we have used $(A\tilde{a}(t)w,w)=t(\tilde{a}(t)w,w)$ (which
follows from $(tI-A)\tilde{a}(t)=q_{A}(t)I$ and $(w,w)=0$) and
\begin{equation}\label{patrat}
\frac{(\tilde{a}(t)w,\tilde{a}(t)w)}{(Aw,w)}=
q_{A}^{\prime}(t)p_{A,x}(t)-q_{A}(t)p_{A,x}^{\prime}(t),
\end{equation}
which has been proved in \cite{paul}.

\end{proof}

\section{Definition of generalised K\"{a}hler cones; the
Bochner-flatness condition}\label{det}

Let $(N,H,I)$ be an oriented CR manifold and $\{
T_{r}\}_{r\in\mathcal J}$, with $\mathcal
J\subseteq\mathbb{R}^{>0}$ a connected open interval, a family of
Reeb vector fields of Sasaki structures on $(N,H,I)$, with contact
forms $\{\theta\}_{r\in \mathcal J}.$ Let
$\omega_{r}:=\frac{1}{2}d\theta_{r}\in\Lambda^{2}(N)$ and
$g_{r}:=\omega_{r}(\cdot ,I\cdot )$ be the corresponding (positive
definite) metrics of the contact bundle $H$. On the cone manifold
$N\times\mathcal J$ define the vector fields $T$, $V$, a complex
structure $J$ and a $2$-form $\omega$ as in the Introduction, and
also a symmetric bilinear form $g:=\omega (\cdot ,J\cdot ).$

Let $M\subset N\times\mathcal J$ be the open subset where $g(T,T)$
is positive. Define a positive function $f:M\to \mathbb{R}^{>0}$
by $g(T,T)=r^{2}f$. Note that the restriction $f_{r}:=f(\cdot ,r)$
of $f$ to $N_{r}:=M\cap N\times \{ r\}$ is positive. We introduce
a new family of contact forms
$\tilde{\theta}_{r}=\frac{1}{f_{r}}\theta_{r}$; for any $r$, the
contact form $\tilde{\theta}_{r}$ is defined on $N_{r}$ (viewed as
an open subset of $N$).\\

\textbf{Conventions of notations:}
For a function $h:M\to \mathbb{R}$, we shall denote 
by $\dot{h}:M\to\mathbb{R}$ its derivative with respect to $r$ and 
by $h_{r}:N_{r}\to\mathbb{R}$ the restriction of $h$ 
to a level set $N_{r}.$

\begin{lem}\label{functionG}
The following equalities hold:
\begin{equation}\label{functf}
f=1+\frac{r\dot{\theta}_{r}(T_{r})}{2}
\end{equation}
and
\begin{equation}\label{tilde}
\dot{\tilde{\theta}}_{r}=-\frac{2G}{r}\tilde{\theta}_{r},
\end{equation}
where
$$
G:M\to \mathbb{R},\quad G:=\frac{r\dot{f}}{2f}-f+1.
$$
\end{lem}

\begin{proof}
Note that $\omega =\frac{1}{2}(dr\land d^{J}r+rdd^{J}r).$ It is
straightforward to see that $d^{J}r=r\theta$, where
$\theta\in\Lambda^{1}\left( N\times\mathcal J\right)$ is defined
by $\theta (Z):=\theta_{r}(\pi_{*}Z)$, for a tangent vector $Z\in
T_{(p,r)}\left( N\times\mathcal J\right)$, where $\pi
:N\times\mathcal J\to N$ is the natural projection. At a tangent
space $T_{(p,r)}\left( N\times\mathcal J\right)
=T_{p}N\times\mathbb{R}$,
$$
d\theta =d\theta_{r}+dr\land \dot{\theta}_{r}
$$
and then, restricted to the same tangent space,
\begin{equation}\label{functg}
\omega=r\left( dr\land
\left(\theta_{r}+\frac{r}{2}\dot{\theta}_{r}\right)
+\frac{r}{2}d\theta_{r}\right) .
\end{equation}
It follows that $$ r^{2}f=\omega (V,T)=r^{2}\left(
1+\frac{r\dot{\theta}_{r}(T_{r})}{2}\right) ,
$$
which implies (\ref{functf}). To prove (\ref{tilde}), we take the
derivative with respect to $r$ of the equality
$\tilde{\theta}_{r}=\frac{1}{f_{r}}\theta_{r}$ and we use the fact
that $\dot{\theta}_{r}=\dot{\theta}_{r}(T_{r})\theta_{r}.$ We get
$$
\dot{\tilde{\theta}}_{r}=-\frac{\dot{f}_{r}}{f_{r}}\tilde{\theta}_{r}+\frac{1}{f_{r}}
\dot{\theta}_{r}=\left(
-\frac{\dot{f}_{r}}{f_{r}}+\dot{\theta}_{r}(T_{r})\right)\tilde{\theta}_{r}
=-\frac{2G}{r}\tilde{\theta}_{r},
$$
which proves our Lemma.
\end{proof}

\begin{lem}\label{kahler}
The pair $(\omega ,J)$ defines a K\"{a}hler structure on
$M$.\end{lem}

\begin{proof}
From relation (\ref{inte}), it is clear that the integrability
tensor $N^{J}$ of the almost complex structure $J$, applied to a
pair of sections $(X,Y)$ of $H$, vanishes. On the other hand,
$N^{J}(X,V)$, restricted to a level set $N\times \{r\}$, is equal
to $-L_{{T}_{r}}(I)(X)$, which is zero, because $T_{r}$ is the
Reeb vector field of a Sasaki structure. It follows that $J$ is
integrable. From (\ref{functg}) it is easy to see that $T$ is
hermitian orthogonal to $H$ and that the restriction of $g$ to
$H\subset T_{(p,r)}(N\times{\mathcal J})$, coincides with
$r^{2}g_{r}$, which is positive definite. We deduce that $g$ is
positive definite on the subset $M$ of $N\times\mathcal J$, 
where $g(T,T)>0$, and that $(M,\omega ,J)$ is
a K\"{a}hler manifold (the $2$-form $\omega$ being closed).
\end{proof}

\begin{defn} The K\"{a}hler manifold
$(M,\omega ,J)$ is a generalised K\"{a}hler cone
over the CR manifold $(N,H,I)$.
It is a restricted generalised K\"{a}hler cone if the function
$f$ is constant along the trajectories of the vector field $T$.
\end{defn}

\textbf{Convention:} For simplicity, in this paper we will
consider only restricted generalised K\"{a}hler cones; when we
refer to a generalised K\"{a}hler cone, we will actually mean
restricted generalised K\"{a}hler cone; this is true also for the
statement of Theorem \ref{main}.

\begin{rem}{\rm {\bf Main class of generalised K\"{a}hler cones:}
We shall be mainly concerned with generalised K\"{a}hler cones
over (open subsets) of hermitian CR spheres. Suppose that
$N\subset\Sigma^{2m+1}$ is an open subset of the hermitian CR
sphere of complex null lines in $W=\mathbb{C}^{m+1,1}$. Then
$\theta_{r}=\frac{\eta}{(B_{r}w,w)}$,
$\tilde{\theta}_{r}=\frac{\eta}{(A_{r}w,w)}$ for some hermitian
trace-free operators $A_{r}$, $B_{r}$ of $W$. The condition
$T(f)=0$ is equivalent to $[A_{r},B_{r}]=0$ for any $r$, as the
following Lemma shows:}
\end{rem}

\begin{lem}\label{operatori}

\begin{enumerate}

\item The operators $A_{r}$ and $B_{r}$
are related in the following way:
\begin{equation}\label{difbr}
A_{r}=B_{r}-\frac{r}{2}\dot{B}_{r}.
\end{equation}

\item The functions $f$ and $G$ have the following expression:
for any $(x,r)\in M$,
$$
f_{r}(x)=\frac{(A_{r}w,w)}{(B_{r}w,w)},\quad
G(x,r)=\frac{r(\dot{A}_{r}w,w)}{2(A_{r}w,w)};\quad w\in x,\quad
w\neq 0.
$$

\item The condition $T(f)=0$ is equivalent to $[A_{r},B_{r}]=0$
for any $r$.

\end{enumerate}

\end{lem}

\begin{proof}
Note that
$\dot{\theta}_{r}=-\frac{(\dot{B}_{r}w,w)}{(B_{r}w,w)}\theta_{r}$
and
$\dot{\theta}_{r}(T_{r})=-\frac{(\dot{B}_{r}w,w)}{(B_{r}w,w)}.$
It follows that
\begin{equation}\label{f}
f=1+\frac{r\dot{\theta}_{r}(T_{r})}{2}=
1-\frac{r}{2}\frac{(\dot{B}_{r}w,w)}{(B_{r}w,w)}
=\frac{\left( B_{r}w-\frac{r}{2}\dot{B}_{r}w,w\right)}{(B_{r}w,w)}.
\end{equation}
We deduce that
\begin{equation}\label{uitat}
\tilde{\theta}_{r}=\frac{1}{f_{r}}\theta_{r}=
\frac{(B_{r}w,w)}{\left( B_{r}w-\frac{r}{2}
\dot{B}_{r}w,w\right)}\cdot \frac{\eta}{(B_{r}w,w)}
=\frac{\eta}{(B_{r}w-\frac{r}{2}\dot{B}_{r}w,w)}.
\end{equation}
Relation (\ref{difbr}) follows from (\ref{uitat}),
$\tilde{\theta}_{r}=\frac{\eta}{(A_{r}w,w)}$ and Lemma
\ref{identitate}. From (\ref{difbr}) and (\ref{f}) we get the
expression of $f.$ On the other hand, from Lemma \ref{functionG}
and $\tilde{\theta}_{r}=\frac{\eta}{(A_{r}w,w)}$, we have
$$
\dot{\tilde{\theta}}_{r}=-\frac{2G}{r}
\tilde{\theta}_{r}=
-\frac{(\dot{A}_{r}w,w)}{(A_{r}w,w)}\tilde{\theta}_{r},$$
which implies that $G$ is of the required form. To prove the last statement,
note that
\begin{align*}
T_{r}(f_{r})&=\frac{[(A_{r}T_{r}w,w)+
(A_{r}w,T_{r}w)]
(B_{r}w,w)-[(B_{r}T_{r}w,w)+(B_{r}w,T_{r}w)]
(A_{r}w,w)}{(B_{r}w,w)^{2}}\\
&=\frac{i([A_{r},B_{r}]w,w)}{(B_{r}w,w )},
\end{align*}
since, at a point $x\in N$,
$T_{r}(x)\in\mathrm{Hom}(x,x^{\perp}/x)$ is the homomorphism
$T_{r}w=iB_{r}w$ $\mathrm{mod} w$ and the operators $A_{r}$ and
$B_{r}$ are hermitian. We conclude from Lemma \ref{identitate}.

\end{proof}

\begin{lem}\label{LC}
The Levi-Civita connection $D^{g}$ of a generalised K\"{a}hler cone
$(M,\omega ,J)$
has the following expression:

\begin{align*}
D^{g}_{X}Y&=\nabla^{r}_{X}Y-
\omega_{r} (X,Y)T_{r}-g_{r}(X,Y)V\\
D^{g}_{V}Y&=fY+\frac{Y(f)}{2f}V-\frac{(JY)(f)}{2f}T\\
D^{g}_{T}Y&={\mathcal L}_{T}(Y)+fJY+\frac{Y(f)}{2f}T+\frac{(JY)(f)}{2f}V\\
D^{g}_{X}V&=fX+\frac{X(f)}{2f}
V-\frac{(JX)(f)}{2f}T\\
D^{g}_{V}V&=-\frac{1}{2}v+( G+f)V\\
D^{g}_{T}V&=\frac{1}{2}Jv+(G+3f-2)T.
\end{align*}
Here $X,Y\in\Gamma (H)$, the vector field $D^{g}_{X}Y$ is
restricted to a level set $N_{r}$, $\nabla^{r}$ is the Tanaka
connection of the contact form $\theta_{r}$ of the CR manifold
$(N,H,I)$ and $v$ is a vector field on $M$ which belongs, at any
point $(p,r)\in M$, to $H_{p}\subset T_{(p,r)}M$ and is determined
by the condition:
\begin{equation}\label{v0}
g_{r}(X,v_{(p,r)})=df(X),\quad \forall
X\in H_{p}\subset T_{(p,r)}M.
\end{equation}

\end{lem}

\begin{proof}
The proof is a straightforward computation based
on the Koszul formula. It uses the expression of the Tanaka connection
on the contact bundle $H$, mentioned in Section \ref{basic}.
\end{proof}

\begin{lem}\label{curb}
The curvature $R^{g}$ of a generalised K\"{a}hler cone
$(M,\omega ,J)$ has
the following expression:
\begin{align*}
g(R^{g}_{X,T}Y,Z)&=-\frac{Y(f)}{2}\omega
(X,Z)-\frac{(JY)(f)}{2}g(X,Z)+
\frac{Z(f)}{2}\omega (X,Y)\\
&+\frac{(JZ)(f)}{2}g(X,Y)-X(f)\omega (Y,Z)\\
g(R^{g}_{X_{1},X_{2}}Y,Z)&=g(R^{\nabla^{r}}_{X_{1},X_{2}}Y,Z)-
\frac{f}{r^{2}}\{ g(X_{1},Y)g(X_{2},Z)-g(X_{2},Y)g(X_{1},Z)\}\\
&+\frac{f}{r^{2}}\{-\omega (X_{1},Y) \omega (X_{2},Z) + \omega
(X_{2},Y)\omega (X_{1},Z)-
2\omega (X_{1},X_{2})\omega (Y,Z)\}\\
g(R^{g}_{X,T}V,Y)&=\frac{{r}^{2}}{2}( \nabla^{r}df)^{J,-}(X,JY)+\frac{r^{2}f}{2}
\left(\nabla^{r}\left(\frac{df}{f}\right)\right)^{J,+}(X,JY)\\
&+f(2-G-2f)\omega (X,Y)\\
g(R^{g}_{T,V}V,T)&=g(v,v)+r^{2}f\left( (G-2)(8f-2)+r\dot{G}+12f^{2}\right)\\
g(R^{g}_{T,V}V,Z)&=-\frac{r^{3}}{2}d\dot{f}(JZ)+
r^{2}(G-1)df(JZ),
\end{align*}
where $X,X_{1},X_{2},Y,Z\in\Gamma (H)$,
$g(R^{g}_{X_{1},X_{2}}Y,Z)$ is restricted to a level set $N_{r}$,
$v$ is the vector field defined by (\ref{v0}) and the superscripts
$J,+$ and $J,-$ denote the $J$-invariant part, respectively the
$J$-anti-invariant part of a bilinear form.

\end{lem}

\begin{proof}
The proof is a lenghty but straightforward computation.
\end{proof}

\begin{prop}\label{flat}
The generalised K\"{a}hler cone $(M,\omega ,J)$ is Bochner-flat if and only if
the CR manifold $(N,H,I)$ is flat and the following
two conditions hold:

\begin{enumerate}

\item The function $G$ depends only on $r$.

\item For every $r\in\mathcal J$, the
contact form $\tilde{\theta}_{r}
=\frac{1}{f_{r}}{\theta}_{r}$ is the contact form of a Sasaki
structure on $(N_{r},H,I)$, which determines an
Einstein K\"{a}hler structure on the quotient $N_{r}/\tilde{T}_{r}$
(where $\tilde{T}_{r}$ is the Reeb vector field of $\tilde{\theta}_{r}$),
with modified Ricci tensor
$$
\tilde{S}_{r}=\left(
\frac{r\dot{G}}{2}-G+2\right)\mathrm{Id}.
$$
\end{enumerate}

\end{prop}

\begin{proof}
The K\"{a}hler manifold $(M,\omega ,J)$  is Bochner flat if and
only if
\begin{equation}\label{cond}
R^{g}=c_{\mathcal K}^{*}(S),
\end{equation}
for a tensor field $S\in\mathrm{Sym}^{1,1}(M)$. Plugging into
(\ref{cond}) the arguments $(T,V,V,T)$ and $(T,V,V,Z)$, for
$Z\in\Gamma (H)$, and using formula (\ref{dual}) for the adjoint
of the Ricci contraction, we readily deduce that $S(T,T)$ and
$S(Z,T)$ are related to the curvature $R^{g}$ as follows:
\begin{equation}\label{conditie}
S(T,T)=-\frac{1}{2r^{2}f}g(R^{g}_{T,V}V,T);\quad
S(Z,T)=-\frac{1}{r^{2}f} g(R^{g}_{T,V}V,Z).
\end{equation}
On the other hand, from Lemma \ref{curb} we know that
\begin{align*}
g(R^{g}_{X,V}Y,Z)&=\frac{Y(f)g(X,Z)}{2}-\frac{Z(f)g(X,Y)}{2}
-\frac{(JY)(f)\omega (X,Z)}{2}\\
&+\frac{(JZ)(f)\omega (X,Y)}{2}-(JX)(f)\omega (Y,Z),
\end{align*}
for every $X,Y,Z\in\Gamma (H)$,
and, since $c_{\mathcal K}^{*}(S)(X,V,Y,Z)
=g(R^{g}_{X,V}Y,Z)$, we readily get
$S(Z,T)=2(JZ)(f).$
Combining this with the second relation (\ref{conditie}), we deduce that
$$
-\frac{1}{r^{2}f}g(R^{g}_{T,V}V,Z)=2(JZ)(f),
$$
which is equivalent, using the expression of $g(R^{g}_{T,V}V,Z)$
provided by Lemma \ref{curb}, to $(JZ)(G)=0.$ Since $Z\in\Gamma
(H)$ is arbitrary, we obtain the first condition of the
Proposition (since $X(G)=Y(G)=0$, for $X,Y\in\Gamma (H)$, also
$[X,Y](G)=0$; recall now that vector fields of the form
$\{X,[X,Y],X,Y\in\Gamma (H)\}$ span the entire $TN$). To obtain
the second condition of the Proposition, we notice that the
expression of $g(R^{g}_{X_{1},X_{2}}Y,Z)$ found in Lemma
\ref{curb}, combined with (\ref{cond}), imply that the CR manifold
$(N,H,I)$ is flat and that on the bundle $H$ restricted to a level
set $N_{r}$,
\begin{equation}\label{doi}
S^{\perp}=\frac{1}{r^{2}}\left( S_{r}-2f\mathrm{Id}\right)
\end{equation}
where $S^{\perp}:H\to H$ is induced by $S$ followed by
$g$-orthogonal projection and
$S_{r}\in\mathrm{End}(H)$ is the modified Ricci tensor of the
K\"{a}hler curvature
$R^{\nabla^{r}}\in\Lambda^{1,1}(H)\odot\Lambda^{1,1}(H)$ of
$\nabla^{r}$.
Plugging into (\ref{cond}) the argument $(X,T,V,Y)$
and using relation
(\ref{doi}), we obtain
\begin{align*}
g(R^{g}_{X,T}V,Y)&=\frac{1}{4}\left( S(JX,Y)g(V,V)
+\omega (X,Y)S(V,V)\right)\\
&=\frac{f}{4}g(S_{r}(JX),Y)-\frac{1}{4}\left( 2f^{2}
+\frac{g(R^{g}_{T,V}V,T)}{2r^{2}f}\right) \omega (X,Y).
\end{align*}
Using the expression of $g(R^{g}_{X,T}V,Y)$ provided by
Lemma \ref{curb} we deduce that
\begin{align*}
g_{r}(S_{r}X,Y)&=-\frac{2}{f}\nabla^{r}(df)^{J,-}(X,Y)-
2\nabla^{r}\left(\frac{df}{f}\right)^{J,+}(X,Y)\\
&+\left( 4(2-G)-6f+\frac{g(R^{g}_{T,V}V,T)}{2f^{2}r^{2}}\right) g_{r}(X,Y).
\end{align*}
This relation clearly implies that $\nabla^{r}(df)^{J,-}\vert_{H\times H}=0$,
which means that $\tilde{\theta}_{r}=\frac{1}{f_{r}}\theta_{r}$
(for any $r$) is the contact form of a Sasaki structure \cite{pl}.
Moreover, the modified Ricci tensor $\tilde{S}_{r}$ of the
the Sasaki structure determined by $\tilde{\theta}_{r}$ is related to $S_{r}$
in the following way \cite{liana}
\begin{equation}
\frac{1}{f}g_{r}(\tilde{S}_{r}(X),Y)=g_{r}(S_{r}(X),Y)
+2\nabla^{r}\left( \frac{df}{f}\right)^{J,+}
(X,Y)-\frac{g_{r}(v,v)}{2f^{2}}g_{r}(X,Y).
\end{equation}
We deduce, using the previous expression of 
$g_{r}(S_{r}(X),Y)$, that
$$
\frac{1}{f}g_{r}(\tilde{S}_{r}(X),Y)+
\left(\frac{g_{r}(v,v)}{2f^{2}}
+4(G-2)+6f
-\frac{g(R^{g}_{T,V}V,T)}{2r^{2}f^{2}}\right) g_{r}(X,Y)=0.
$$
Using again Lemma \ref{curb} for the expression of
$g(R^{g}_{T,V}V,T)$ we obtain the second condition of the
Proposition. Conversely, it is easy to check that the two
conditions of the Proposition ensure the Bochner-flatness of
$(M,\omega ,J).$

\end{proof}

The main result of this Section is the following:

\begin{prop}\label{explicit}

Let $(M,g,J)$ be a Bochner-flat generalised K\"{a}hler cone of
complex dimension $m+1\geq 3$, defined by a family of Sasaki Reeb
vector fields $\{ T_{r}\}$ over a CR manifold $(N,H,I)$. Then
$(N,H,I)$ is locally isomorphic to the standard CR sphere
$\Sigma^{2m+1}$ of complex null lines in a complex hermitian
vector space $W$ of signature $(m+1,1)$, and $\{ T_{r}\}$ is
defined by one of the following families of hermitian operators
$B_{r}$ of $W$:

\begin{enumerate}

\item $B_{r}=r^{2}(B-\mu (r^{2})A).$ Here the real function $\mu$ satisfies
${\mu}^{\prime}>0$ and is a solution
of the differential equation
\begin{equation}\label{muec}
{\mu}^{\prime}=\frac{1}{2}\mu^{2}+d,
\end{equation}
where $d\in\mathbb{R}$ is an arbitrary real number. The operator
$A$ is a hermitian semi-simple operator, with a positive definite
eigenspace, of dimension $m+1$, which corresponds to the
eigenvalue $\frac{1}{2(m+2)}$ and a $1$-dimensional timelike
eigenspace, which corresponds to the eigenvalue
$-\frac{(m+1)}{2(m+2)}$.

\item $B_{r}=r^{2}(B+ \mu (r^{2})A)$, where $\mu$ satisfies  (\ref{muec}) and
${\mu}^{\prime}<0.$ The operator $A$ is semi-simple, with an eigenspace
of signature $(m,1)$, which corresponds to the eigenvalue
$-\frac{1}{2(m+2)}$, and a $1$-dimensional spacelike eigenspace,
which corresponds to the eigenvalue $\frac{m+1}{2(m+2)}.$

\item $B_{r}=r^{2}(B-r^{2}A)$, where $A$ is $1$-step parabolic, with all
eigenvalues equal to zero.

\item $B_{r}=r^{2}\left(B -\frac{e^{\lambda r^{2}}}{\lambda} A\right)$, where
$\lambda\in\mathbb{R}\setminus \{ 0\}$, $A$ is $1$-step parabolic
with all eigenvalues equal to zero.

\end{enumerate}

In all these cases, $B$ is any hermitian, trace-free operator
of $W$ which commutes with $A$ (see Remark \ref{op}).

\end{prop}

\begin{proof}
Since $(N,H,I)$ is flat (see Proposition \ref{flat}) and of
dimension bigger than three, we can assume, restricting $N$ if
necessary, that $(N,H,I)$ is an open subset of the hermitian CR
sphere $\Sigma^{2m+1}$ of complex null lines in a complex
hermitian vector space $W$ of signature $(m+1,1)$ \cite{vesti}. As
explained in Lemma \ref{operatori}, the two families of contact
forms $\{\theta_{r}\}_{r\in\mathcal J}$ and
$\{\tilde{\theta}_{r}\}_{r\in\mathcal J}$ are generated by two
families of hermitian trace-free endomorphism $\{
B_{r}\}_{r\in\mathcal J}$ and $\{ A_{r}\}_{r\in\mathcal J}$ of $W$
respectively such that, for any $r$, the operators $A_{r}$ and
$B_{r}$ commute (see Lemma \ref{operatori}). From Lemma
\ref{functionG} we know that
$\dot{\tilde{\theta}}_{r}=-\frac{2G_{r}}{r}\tilde{\theta}_{r}$.
Since $G$ depends only on $r$,  we get $\tilde{\theta}_{r}
=e^{-\int_{r_{0}}^{r}\frac{2G_{q}}{q}dq}\tilde{\theta}_{r_{0}}$
and we infer that the modified Ricci tensor $\tilde{S}_{r}$ of
$\tilde{\theta}_{r}$ has the expression
$\tilde{S}_{r}=e^{\int_{r_{0}}^{r}\frac{2G_{q}}{q}dq}\tilde{S}_{r_{0}}$.
The second condition of Theorem \ref{flat} is equivalent with
$$
(\frac{r\dot{G}}{2}-G+2)
\mathrm{Id}=e^{\int_{r_{0}}^{r}\frac{2G_{s}}{s}ds}
\tilde{S}_{r_{0}}
$$
and implies
\begin{equation}\label{ecG}
\left(\frac{r\dot{G}}{2}-G+2\right)^{\prime}
=\frac{2G}{r}\left(\frac{r\dot{G}}{2}-G+2\right) .
\end{equation}
Equation (\ref{ecG}) can be solved as follows: define
a real function $\mu$ in the following way:
\begin{equation}\label{mu}
\mu (t)=\frac{G(\sqrt{t})-2}{t}.
\end{equation}
We shall write equation (\ref{ecG}) in terms of $\mu .$ For this,
we first take the derivative of $r^{2}\mu (r^{2})=G(r)-2$ and we
get:
$$
\dot{G}(r)=2r\mu (r^{2})+2r^{3}{\mu}^{\prime}(r^{2}).
$$
It easily follows that
\begin{equation}\label{G1}
\frac{r\dot{G}}{2}-G+2=r^{4}{\mu}^{\prime}(r^{2}).
\end{equation}
Equation (\ref{ecG}) becomes
${\mu}^{\prime\prime}={\mu}^{\prime}\mu$. Since $\mathcal J$ is connected,
$\mu$ satisfies (\ref{muec}), for a constant $d\in\mathbb{R}.$
We have the following three possibilities:

\begin{enumerate}

\item  ${\mu}^{\prime}>0.$ From Lemma \ref{operatori} we deduce that
$$
(\dot{A}_{r}w,w)=\frac{2}{r}G(r)(A_{r}w,w)=
\frac{2}{r}\left( r^{2}\mu (r^{2})+2\right) (A_{r}w,w).
$$
Since $\int r\mu
(r^{2})dr=\frac{1}{2}\mathrm{ln}\left({\mu}^{\prime}(r^{2})\right)$
when ${\mu}^{\prime}>0$ we get
$$
(A_{r}w,w)=K(w)r^{4}{\mu}^{\prime}(r^{2}),
$$
where $K=K(w)$ depends only on $w$.
Equivalently,
\begin{equation}\label{11}
A_{r}={\mu}^{\prime}(r^{2})r^{4}A,
\end{equation}
where $A\in\mathrm{End}(W)$ is hermitian, trace-free, satisfies
$(Aw,w)>0$, for $w\in x$ non-zero, when $(x,r)\in M$. Moreover,
(\ref{G1}) together with the second condition of Proposition
\ref{flat} imply that the modified Ricci tensor $S_{A}$ of the
K\"{a}hler-Einstein structure $M_{A}$ is the identity
endomorphism, from where we deduce that $A$ is as in the statement
of the Proposition (see \cite{paul}). On the other hand, from
Lemma \ref{operatori}, $B_{r}$ must satisfy (\ref{difbr}), with
$A_{r}={\mu}^{\prime}(r^{2})r^{4}A$. It follows that
$$
B_{r}=r^{2}\left( B - \mu (r^{2})A\right) ,
$$
where $B\in\mathrm{End}(W)$ is hermitian and trace free.\\

\item ${\mu}^{\prime}<0.$
Then $\int r\mu (r^{2})dr =\frac{1}{2}\mathrm{ln}\left( -{\mu}^{\prime}(r^{2})\right).$
A similar argument shows that
\begin{align*}
A_{r}&=-r^{4}{\mu}^{\prime}(r^{2})A\\
B_{r}&=r^{2}\left( B+\mu (r^{2})A\right),
\end{align*}
but in this case the Bochner-flat K\"{a}hler structure $M_{A}$ has
the modified Ricci operator $S_{A}=-\mathrm{Id}$, which implies
that $A$ is as in the statement of the Proposition (see
\cite{paul}).

\item It remains to consider the case when
the function $\mu$ is constant. Then $\mu (t)=\lambda$ for
$\lambda\in\mathbb{R}$, $G(x,r)=\lambda r^{2}+2$ and
$$
(\dot{A}_{r}w,w)= \frac{2G}{r}(A_{r}w,w)= \frac{2(\lambda
r^{2}+2)}{r}(A_{r}w,w).
$$
We distinguished two subcases: (i) $\lambda =0$; (ii) $\lambda\neq 0.$
In subcase (i) we obtain
$$
A_{r}=r^{4}A,\quad B_{r}=r^{2}(B-r^{2}A),
$$
and in subcase (ii),
$$
A_{r}=r^{4}e^{\lambda r^{2}}A,\quad
B_{r}=r^{2}\left( B-\frac{e^{\lambda r^{2}}}{\lambda} A\right) .
$$
Since $\frac{r\dot{G}}{2}-G+2=0$,
the K\"{a}hler structure $M_{A}$ is flat and hence
the endomorphism $A$ is $1$-step parabolic, with all eigenvalues
zero (see \cite{paul}).

\end{enumerate}

In all cases (1), (2) and (3), the generalised K\"{a}hler cone
condition $T(f)=0$ becomes $[A,B]=0$ (see Lemma \ref{operatori}).

\end{proof}

\begin{rem}\label{op}{\rm
The condition $[A,B]=0$ of Proposition \ref{explicit} determines
the operator $B$ as follows:

\begin{enumerate}

\item In the first case of Proposition \ref{explicit}, $B$ preserves,
up to a multiplicative constant, a timelike eigenvector $v$
(unique, up to a non-zero multiplicative constant) of $A$. On the hermitian
orthogonal $v^{\perp}$, the hermitian metric $(\cdot ,\cdot )$ is
positive definite, $A:v^{\perp} \to v^{\perp}$ is a multiple of
the identity endomorphism and $B :v^{\perp}\to v^{\perp}$, being
hermitian, is diagonalisable. It follows that $A$ and $B$ are
simultaneously  diagonalisable.

\item In the second case of Proposition \ref{explicit}, $B$ preserves,
up to a  multiplicative constant, a spacelike eigenvector $v$
(unique, up to a non-zero multiplicative constant) of $A$, which
corresponds to the eigenvalue $\frac{m+1}{2(m+2)}.$ On the
hermitian orthogonal $v^{\perp}$, the hermitian metric $(\cdot
,\cdot )$ has signature $(m,1)$, $A$ is a multiple of the identity
endomorphism and $B:v^{\perp}\to v^{\perp}$ can be elliptic,
hyperbolic, $1$- or $2$-step parabolic.

\item Consider now the cases three and four of Proposition \ref{explicit}.
Note that $A=0$ on any positive definite eigenspace of $B$ (since
$[A,B]=0$, $A$ preserves such an eigenspace, say $W_{1}$, of $B$;
because the hermitian metric $(\cdot ,\cdot )$ is positive
definite on $W_{1}$ and $A$ is hermitian, $A$ is diagonalisable on
$W_{1}$; this forces it to be zero, because $A$ does not have
non-zero eigenvalues). Let us denote by $W_{1},\cdots ,W_{s}$, the
positive definite eigenspaces of $B$ and by $W_{0}$ the hermitian
orthogonal of the direct sum $\oplus_{j=1}^{s}W_{j}$. The
eigenspaces $W_{j}$ (for $j\in \{ 1,\cdots ,s\}$) correspond to
eigenvalues, say $\beta_{j}$, of $B$, which can be any real
numbers. It remains to study the restriction $B_{0}$ of $B$ to
$W_{0}$. We notice first that $B_{0}$ cannot be hyperbolic: if it
was hyperbolic, then $B$ would have two complex non-real
eigenvalues, say $\delta$ and $\bar{\delta}$, with $1$-dimensional
eigenspaces, generated by two null independent vectors $v_{1}$ and
$v_{2}$ respectively. However, since $[A,B]=0$, $BAv_{1}=\delta
Av_{1}$ and $BAv_{2}=\bar{\delta}Av_{2}$, which imply that
$Av_{1}=Av_{2}=0$ (because $A$ has no non-zero eigenvalues). But
if we take an orthonormal basis $\{ e_{0},\cdots ,e_{n}\}$ of
$W_{0}$ in which
$$
A=\left(\begin{tabular}{ccccc}
$-1$ & $1$ & $0$ & $\cdots$ & $0$\\
$-1$ & $1$ & $0$ & $\cdots$ & $0$\\
$0$ & $0$ & $0$ & $\cdots$ & $0$\\
$\cdots$ & $\cdots$ & $\cdots$ & $\cdots$ & $\cdots$\\
$0$ & $0$ & $0$ & $0$ & $0$
\end{tabular}\right) ,
$$
(which is possible since $A:W_{0}\to W_{0}$ is $1$-step parabolic with 
all eigenvalues equal to zero)
the conditions $v_{1}, v_{2}$ null and $Av_{1}=Av_{2}=0$ would
imply that $v_{1}$ and $v_{2}$ are multiples of $e_{0}+e_{1}$. In
particular they would be dependent, which is a contradiction. We
conclude that $B_{0}$ can be elliptic or $1$- or $2$-step
parabolic. Therefore, $B_{0}=\gamma I+N$, for an endomorphism $N$
of $W_{0}$ which commutes with $A$ and which satisfies $N^{3}=0$,
and for $\gamma\in\mathbb{R}$ which is different from all
$\beta_{j}$. The endomorphism $N^{\perp}$ of
$\widehat{{W}_{0}}:=\mathrm{Span}\{ e_{2},\cdots ,e_{n}\}$
obtained from $N$ by restriction and orthogonal projection, is
hermitian on $\widehat{{W}_{0}}$. Because the metric $(\cdot
,\cdot )$ is positive definite on $\widehat{{W}_{0}}$, $N^{\perp}$
is diagonalisable and hence there is a basis $\{
e_{2}^{\prime},\cdots ,e_{n}^{\prime}\}$ of $\widehat{{W}_{0}}$,
such that $N^{\perp}$ is diagonal in this basis. If we consider
now the basis ${\mathcal B}: =\{ e_{0},e_{1},e_{2}^{\prime},\cdots
,e_{n}^{\prime}\}$ of $W_{0}$, it is straightforward to see that
$[A,N]=0$ and $N$-hermitian imply that
$$
N=\left(\begin{tabular}{cccccc}
$\gamma _{0}$ & $\alpha$ & $\mu_{2}$ & $\mu_{3}$ & $\cdots$ & $\mu_{n}$\\
$-{\alpha}$ & $\gamma_{1}$ & $\mu_{2}$ & $\mu_{3}$ & $\cdots$ & $\mu_{n}$\\
$-\bar{\mu}_{2}$ & $\bar{\mu}_{2}$ & $\gamma _{2}$ & $0$ & $\cdots$ & $0$\\
$-\bar{\mu}_{3}$ & $\bar{\mu}_{3}$  & $0$ & $\gamma _{3}$ & $\cdots$ & $0$\\
$\vdots$ & $\vdots$ & $\vdots$ & $\vdots$ & $\vdots$ & $\vdots$\\
$-\bar{\mu}_{n}$ & $\bar{\mu}_{n}$ & $0$ & $0$ & $\cdots$ & $\gamma _{n}.$
\end{tabular}\right)
$$
in the basis $\mathcal B .$ Moreover, $N^{3}=0$ if and only if
$\gamma_{0}=-\gamma_{1}=- \alpha$ and $\gamma_{j}=0$ for any $j\in
\{ 2,\cdots ,n\}$ and $N^{2}=0$ if and only
$\gamma_{0}=-\gamma_{1} = - \alpha $, $\gamma_{j}=\mu_{j}=0$, for
any $j\in \{ 2,\cdots ,n\}$. Since $B$ is trace-free, 
the real constants $\beta_{j}$ and $\gamma$ must satisfy
$(n+1)\gamma
+\sum_{i=1}^{s}n_{i}\beta_{i}=0$, where $n_{i}$ is the dimension
of $W_{i}$ (and $n+1$ is the dimension of $W_{0}$).
\end{enumerate}}

\end{rem}

\section{The local geometry of Bochner-flat
generalised K\"{a}hler cones}\label{quotient}

In this section we prove our main Theorem \ref{main}. We will do
this by analysing the local types of the Bochner-flat generalised
K\"{a}hler cones determined in Proposition \ref{explicit}. The
results we shall obtain in this Section can be summarized by the
following table:
$$
\begin{tabular}{|c|c|c|}
\hline
$ $ & ${\mu }^{\prime}>0$ & ${\mu}^{\prime}=0$\\
\hline
$d>0$ & $\mathrm{(all)\ hyperbolic}$ & $---$\\
\hline
$d=0$ & $\mathrm{(all)\ 1-step\ parabolic}$ & $\lambda =0;
\mathrm{(all)\ 2-step\ parabolic}$\\
\hline $d<0$ & $\mathrm{(all)\ elliptic}$ & $\lambda
=\pm\sqrt{-2d},\ \mathrm{(all)\ 1,2-step\ parabolic,\
elliptic}$\\
\hline
\end{tabular}
$$
In particular, we show that all elliptic, hyperbolic, $1$ and $2$-step
parabolic Bochner-flat K\"{a}hler manifolds are locally
generalised K\"{a}hler cones, which proves our main Theorem
\ref{main}.\\

\textbf{Convention:}  Without further explanations, we employ the
notations of the previous section.

\subsection{The case $1$ of Theorem \ref{explicit}}

In this Subsection we analyze the first column of the above table.

\begin{prop}\label{clasif}
Let $(M,g,J)$ be a Bochner-flat generalised K\"{a}hler cone which
belongs to the first case of Proposition \ref{explicit}. Then
$(M,g,J)$ is:

\begin{enumerate}
\item  of hyperbolic type, if $d>0$.

\item of elliptic type, if $d<0$.

\item of $1$-step parabolic type, if $d=0$.

\end{enumerate}
Conversely, any Bochner-flat K\"{a}hler manifold which is of
elliptic, hyperbolic or $1$-step parabolic type can be locally
realised as a generalised K\"{a}hler cone, which belongs to the
first case of Proposition \ref{explicit}.
\end{prop}

\begin{proof}
In the first case of Proposition \ref{explicit},
$$M=\{ (x,r)\in \Sigma^{2m+1}\times
\mathbb{R}^{>0}:\quad (Bw,w)>\frac{1}{2}\mu (r^{2})(Aw,w),
(Aw,w)>0, \forall w\in x, w\neq 0\} .
$$
From Remark \ref{op}, there is an orthonormal basis ${\mathcal
B}=\{ e_{0},e_{1},\cdots ,e_{m+1}\}$ of $W$ such that both
operators $A$ and $B$ are diagonal in this basis:
\begin{align*}
B&=\mathrm{diag}(-k,k_{1},\cdots ,k_{m+1})\\
A&=\frac{1}{2(m+2)}\mathrm{diag}(-m-1,1,\cdots ,1).
\end{align*}
Here $k_{j}\in\mathbb{R}$, for any $j\in \{ 1,\cdots ,m+1\}$, and
$k=k_{1}+\cdots +k_{m+1}.$ We shall identify $\Sigma^{2m+1}$ with
the unit sphere $S^{2m+1}$ in $\mathrm{Span}\{ e_{1},\cdots
,e_{m+1}\}$ and $S^{2m+1}\times \mathbb{R}^{>0}$ with
$\mathbb{C}^{m+1}\setminus\{ 0\}$, by means of the diffeomorphism
$$
h:S^{2m+1}\times\mathbb{R}^{>0}\to \mathbb{C}^{m+1} \setminus\{
0\},\quad f(z_{1},\cdots ,z_{m+1},r):=(rz_{1},\cdots ,rz_{m+1}).
$$
Read on the image
$$
h(M)=\{ (z_{1},\cdots ,z_{m+1})\in\mathbb{C}^{m+1}\setminus \{
0\},\quad \sum_{j=1}^{m+1}\left( k+k_{j}-\frac{1}{2}\mu
(r^{2})\right) |z_{j}|^{2}>0\} ,
$$
the complex structure $J$, at a point $z\in h(M)$, satisfies
\begin{align*}
J(V)&=r^{2}\sum_{j=1}^{m+1}
\left( k_{j}+k-\frac{1}{2}\mu (r^{2})\right)\left(
x_{j}\frac{\partial}{\partial y_{j}}-y_{j}\frac{\partial}{\partial x_{j}}
\right)\\
J_{\vert z^{\perp}}&=i
\end{align*}
and the K\"{a}hler form $\omega$ is equal to
$\frac{1}{4}dd^{J}r^{2}$. Here $r^{2}=|z_{1}|^{2}+\cdots
+|z_{m+1}|^{2}$ and $z_{j}=x_{j}+iy_{j}$ for any $j\in\{ 1,\cdots
,m+1\} $. For simplicity, we restrict the K\"{a}hler structure
$(\omega ,J)$ to the subset
\begin{equation}\label{d}
D:=\{ z\in\mathbb{C}^{m+1}\setminus \{ 0\} ,\quad \mu
(r^{2})<2(k_{j}+k),\quad j\in\overline{1,m+1}\} 
\end{equation}
of $h(M).$ We shall consider separately the three cases: $d>0$,
$d<0$ and $d=0.$\\

(1) Suppose that $d=\frac{\beta^{2}}{2}$, with $\beta
>0$. Then $\mu (t)=\beta\mathrm{tg}\left(\frac{\beta
t}{2}\right)$. It can be checked that the map
\begin{equation}\label{F}
F(z_{1},\cdots ,z_{m+1}):=(w_{1}=f_{1}(r^{2})z_{1},\cdots ,
w_{m+1}=f_{m+1}(r^{2})z_{m+1}),
\end{equation}
where
\begin{equation}\label{rh}
f_{j}(t)=\left(\frac{\beta}{\sqrt{2}}\right)^{1/2}
\frac{e^{\frac{(k_{j}+k)t}{2}}}{t^{1/2}{\mu}^{\prime}(t)^{1/4}} ,
\end{equation}
is an isomorphism between the K\"{a}hler manifolds
$(D,\frac{1}{4}dd^{J}r^{2},J)$ and
$(F(D),\frac{1}{4}dd^{J_{0}}x,J_{0})$. Here $J_{0}$ is the
standard complex structure of $\mathbb{C}^{m+1}$ and the positive
function $x=x(w_{1},\cdots , w_{m+1})$ is defined by the implicit
equation
\begin{equation}\label{radial}
\sum_{j=1}^{m+1}\frac{|w_{j}|^{2}}{f_{j}^{2}(x)}=x.
\end{equation}
Let $y=y(w_{1},\cdots ,w_{m+1})$ be related to $x$ by the formula
$$
x=\frac{4}{\beta}\mathrm{arctg}[(1+y)^{1/2}]
$$
and notice that
$$
\frac{e^{(k_{j}+k)x}}{\dot{\mu}(x)^{1/2}}
=\frac{\sqrt{2}|y|}{\beta (2+y)}e^{2l_{j}\mathrm{arctg}[(1+y)^{1/2}]}
$$
where $l_{j}=\frac{2}{\beta}(k_{j}+k)$, for $j\in\{ 1,\cdots
,m+1\} .$ For simplicity, we restrict to the set, say
$D^{\prime}\subset F(D)$, where $y>0.$ On this set,
$(\frac{1}{4}dd^{J_{0}}x,J_{0})$ coincides with the K\"{a}hler
structure
$(\frac{1}{\beta}dd^{J_{0}}\mathrm{arctg}[(1+y)^{1/2}],J_{0})$,
where $y$ is defined by the implicit equation
\begin{equation}\label{h}
\sum_{j=1}^{m+1}
\frac{|w_{j}|^{2}}{e^{2l_{j}\mathrm{arctg}[(1+y)^{1/2}]}}=\frac{y}{2+y}.
\end{equation}
The K\"{a}hler structure
$(\frac{1}{\beta}dd^{J_{0}}\mathrm{arctg}[(1+y)^{1/2}],J_{0})$ is
of hyperbolic type (see Section 2 of \cite{paul1}), isomorphic
with $M_{C}$, where $C$ a hyperbolic hermitian operator of
$\mathbb{C}^{m+2,1}$ with characteristic and minimal polynomials
\begin{align*}
Q_{C}&=\left(\left(t+\frac{2(m+2)k}{\beta
(m+3)}\right)^{2}+1\right) \prod_{j=1}^{m+1}
\left( t-\frac{2}{\beta}\left( k_{j}+\frac{k}{m+3}\right)\right)\\
q_{C}&= \left(\left(t+\frac{2(m+2)k}{\beta ( m+3 )
}\right)^{2}+1\right) \prod_{i=1}^{s} \left(
t-\frac{2}{\beta}\left( k_{i}+\frac{k}{m+3}\right)\right).
\end{align*}
Here $i\in \{ 1,\cdots ,s\}$ parametrises the distinct values of
$\{ k_{1},\cdots ,k_{m+1}\} .$ This proves our first claim.

(2) Next, suppose that $d:=-\frac{\beta^{2}}{2}$, where $\beta
>0$. Then $\mu (t)=\frac{\beta (1+e^{t\beta +p})}{1-e^{t\beta
+p}}$, where $p\in\mathbb{R}.$ By means of the function (\ref{F})
with
\begin{equation}\label{2}
f_{j}(t)
=\frac{\beta^{1/2}2^{1/4}e^{\frac{1}{2}(k_{j}+k)(t+\frac{p}{\beta})
+p}}{t^{1/2}{\mu}^{\prime}(t)^{1/4}}, \quad j\in \{ 1,\cdots ,m+1\} ,
\end{equation}
the K\"{a}hler manifolds $(D,\omega ,J)$ and $(F(D),\omega_{0}
:=\frac{1}{4}dd^{J_{0}}x,J_{0})$ are isomorphic, where
$x=x(w_{1},\cdots ,w_{m+1})$ is defined by the implicit equation
(\ref{radial}), with functions $f_{j}$ defined in (\ref{2}).
Define a new function $y=y(w_{1},\cdots ,w_{m+1})$ by
$$
y=e^{\beta x}-e^{-p}.
$$
For simplicity, we restrict to the subset of $F(D)$, say
$D^{\prime}$, where $y>0.$ On $D^{\prime}$, the function $y$
satisfies the implicit equation
$$
\sum_{j=1}^{m+1}\frac{|w_{j}|^{2}}{e^{p}y(e^{p}y+1)^{\beta_{j}}}=1,
$$
with $\beta_{j}:=\frac{1}{\beta}(k_{j}+k-\frac{\beta }{2}).$ In
terms of $y$, $\omega_{0}
=\frac{1}{4\beta}dd^{J_{0}}\mathrm{log}(e^{p}y+1)$. It follows
that $(\omega_{0} ,J_{0})$ is of elliptic type on $D^{\prime}$
(see Section 2.2 of \cite{paul}) and is isomorphic with the
K\"{a}hler structure $M_{C}$, with $C$ a semi-simple hermitian
operator of $\mathbb{C}^{m+2,1}$,  with eigenvalues $-\frac{\beta
}{2}-\frac{(m+2)k}{m+3},k_{1} +\frac{k}{m+3},\cdots ,
k_{m+1}+\frac{k}{m+3},\frac{\beta }{2}-\frac{(m+2)k}{m+3}.$ Our
second claim follows.

(3) Finally, suppose that $d=0$. Then $\mu (t)=-\frac{2}{t+q}$,
for $q\in\mathbb{R}.$ For simplicity, we assume that $q\geq 0.$
The function (\ref{F}), with
$$
f_{j}(t)=\frac{e^{\frac{(k_{j}+k)t}{2}}(t+q)^{1/2}}{t^{1/2}}
$$
defines an isomorphism between the K\"{a}hler manifolds $(D,\omega
,J)$ and $(F(D),\omega_{0}:=\frac{1}{4}dd^{J_{0}}x,J_{0})$ where
$x =x(w_{1},\cdots ,w_{m+1})$ is a positive function defined
implicitly by the equation
$$
\sum_{j=1}^{m+1}\frac{|w_{j}|^{2}}{e^{(k_{j}+k)x}(x+q)}=1.
$$
Note that the function $y:=x+q>0$ satisfies
\begin{equation}
\sum_{j=1}^{m+1}\frac{|w_{j}|^{2}}{e^{(k_{j}+k)(y-q)}}=y.
\end{equation}
Moreover, $\omega_{0}=\frac{1}{4}dd^{J_{0}}y$. The K\"{a}hler
structure $(\frac{1}{4}dd^{J_{0}}y,J_{0})$ is of $1$-step
parabolic type (see Section 3.1 of \cite{paul}). It is isomorphic
to the K\"{a}hler structure $M_{C}$, where $C$ is a $1$-step
parabolic hermitian operator of $\mathbb{C}^{m+2,1}$, with
characteristic and minimal polynomials
\begin{align*}
Q_{C}(t)&=\left( t+\frac{(m+2)k}{m+3}\right)^{2}\prod_{j=1}^{m+1}
\left( t-k_{j}-\frac{k}{m+3}\right)\\
q_{C}(t)&=\left( t+\frac{(m+2)k}{m+3}\right)^{2}\prod_{i=1,
k_{i}\neq -k}^{s}
\left( t-k_{i}-\frac{k}{m+3}\right) .\\
\end{align*}
As before, $i\in \{ 1,\cdots ,s\}$ parametrises the distinct
values of $\{ k_{1},\cdots ,k_{m+1}\} .$ Our third claim follows.

The last statement of the Proposition follows by an
examination of the minimal and characteristic polynomials of the
operators $C$ we found in each of the cases (1), (2) and (3).

\end{proof}

\subsection{The cases $3$ and $4$ of Proposition \ref{explicit}}\label{45}

In this Section we prove the following

\begin{prop}\label{es}
\begin{enumerate}
\item Let $(M,g,J)$ be a Bochner-flat generalised K\"{a}hler cone,
which belongs to the third case of Proposition \ref{explicit}.
Then $(M,g,J)$ is of $2$-step parabolic type.

\item Let $(M,g,J)$
be a Bochner-flat generalised K\"{a}hler cone, which belongs to
the fourth case of Proposition \ref{explicit}. Then $(M,g,J)$ is
of $2$-step parabolic type, except when $\mu_{j}=0$ for any $j\in
\{ 2,\cdots ,n\}.$ In this case it is of $1$-step parabolic type
if $\alpha\neq 0$ and of elliptic type if $\alpha =0.$

\item Any $2$-step parabolic Bochner-flat K\"{a}hler structure can
be locally realised as a Bochner-flat generalised K\"{a}hler cone
which belongs to the third case of Proposition \ref{explicit}, and
also as a Bochner-flat generalised K\"{a}hler cone which belongs
to the fourth case of Proposition \ref{explicit}.

\end{enumerate}
\end{prop}

We divide the proof into several Lemmas and Propositions. Let
$(M,\omega ,J)$ be a Bochner-flat generalised K\"{a}hler cone,
which belongs to the third or to the fourth case of Proposition
\ref{explicit}. We preserve the notations of Proposition
\ref{explicit} and of Remark \ref{op}.

\begin{prop}\label{expresii}
The Bryant modified Ricci operator $\Theta$ of $(M,g,J)$ has the
following properties:
\begin{align*}
\Theta (L_{j})&=\left(\eta_{j}-\frac{c}{m+3}\right) L_{j}-
\frac{q_{r}(\xi_{j})}{r^{4}\left(\eta_{j}-\gamma \right)}V,\quad j\in \{ 1,\cdots ,l\}\\
\Theta (v_{k})&=
\left(\frac{\lambda_{k}}{r^{2}}-\frac{c}{m+3}\right) v_{k},\quad k\in \{ 1,\cdots ,m-l\}\\
\Theta (V)&=\frac{f}{r^{4}}
\sum_{j=1}^{l}\frac{L_{j}}{p^{\prime}_{r}(\xi_{j})
\left(\eta_{j}-\gamma \right)}-
\left(\frac{f}{r^{2}}-\frac{(m+2)c}{m+3}\right) V.
\end{align*}
Here $\xi_{1},\cdots ,\xi_{l}$ (respectively, $\lambda_{1},\cdots
,\lambda_{m-l})$ are the non-constant (respectively, constant)
eigenvalues of the Bryant modified Ricci operator $\Theta_{r}$ of
$M_{B_{r}}$, $L_{1},\cdots ,L_{l}$ are 
vector fields on $M$ which, at
a point $(x,r)\in M$, belong to $H_{x}=\mathrm{Hom}_{\mathbb{C}}(x,W/x)$
and are the homomorphisms
$L_{j}(w)=\tilde{b}_{r}(\xi_{j})w\quad\mathrm{mod}w$, 
$v_{1},\cdots ,v_{m-l}$ are eigenvectors of
$\Theta$ which correspond to the eigenvalues $\lambda_{1},\cdots
,\lambda_{m-l}$, $p_{r}(t)=\prod_{j=1}^{l}(t-\xi_{j})$ is the
non-constant part of the characteristic polynomial of
$\Theta_{r}$, $\eta_{1}:=\frac{\xi_{1}}{r^{2}},\cdots ,
\eta_{l}:=\frac{\xi_{l}}{r^{2}}$, $q_{r}$ is the minimal
polynomial of $B_{r}$ and $c=\gamma -\lambda$.
\end{prop}

\begin{proof}
Recall that $\Theta$ is related to the modified Ricci tensor $S$
from the proof of Proposition \ref{flat} by
$$
\Theta =\frac{1}{4}\left(
S-\frac{\mathrm{trace}_{\mathbb{R}}(S)}{2(m+3)}\mathrm{Id}\right).
$$
We need to determine $S(L_{j})$, $S(v_{k})$, $S(V)$ and 
$\mathrm{trace}_{\mathbb{R}}(S)$.
From the proof of Proposition \ref{flat} we know that for any
$X\in H$,
\begin{equation}\label{s}
S(X)=\frac{1}{r^{2}}\left( S_{r}X-2fX\right) +\frac{2(JX)(f)}{r^{2}f}T
-\frac{2X(f)}{r^{2}f}V.
\end{equation}
It is easy to check the following equalities:
\begin{equation}\label{df}
df(L_{j})=\frac{2fq_{r}(\xi_{j})}{r^{2}(\eta_{j}-\gamma )};\quad
df(JL_{j})=0;\quad df(v_{k})=0;\quad df(Jv_{k})=0
\end{equation}
which imply, using $\Theta_{r}(L_{j})=\xi_{j}L_{j}$,
$\Theta_{r}(v_{k})=\lambda_{k}v_{k}$ and relation (\ref{theta})
applied to $\Theta_{r}$, that
\begin{align*}
S(L_{j})&=\frac{2}{r^{2}}\left(
2\xi_{j}-\frac{(B_{r}^{2}w,w)}{(B_{r}w,w)} -f\right) L_{j}
-\frac{4q_{r}(\xi_{j})}{r^{4}\left(\eta_{j}-\gamma \right)}V\\
S(v_{k})&=\frac{2}{r^{2}}\left(
2\lambda_{k}-\frac{(B_{r}^{2}w,w)}{(B_{r}w,w)} -f\right) v_{k}.
\end{align*}
To evaluate $S(V)$, we write it as a sum of $S^{\perp}(V)$, the
$g$-orthogonal projection of $S(V)$ on $H$, and
$\frac{S(V,V)}{r^{2}f}V.$ It is easy to check, using the fact that
$S$ is hermitian, relation (\ref{grad}) and $L_{j}=
\frac{p^{\prime}_{r}(\xi_{j})}{2}\mathrm{grad}_{g_{r}}(\xi_{j})$
(see Theorem $2$ of \cite{paul}), that
\begin{equation}\label{sp}
S^{\perp}(V)=\frac{4f}{r^{4}}
\sum_{j=1}^{l}\frac{L_{j}}{p_{r}^{\prime}(\xi_{j})(\eta_{j}-\gamma
)}.
\end{equation}
To evaluate $S(V,V) = S(T,T)$, we use the first
equality of (\ref{conditie}) and the expression of
$g(R^{g}_{T,V}V,T)$ provided by Lemma \ref{curb}. Notice that the
vector field $v$ defined in Lemma \ref{LC} has the
following expression: at a point $(x,r)\in M$, $v_{(x,r)}\in H_{x}
=\mathrm{Hom}(x,W/x)$ is equal to 
$$
v_{(x,r)}(w)=2\left( A_{r}w-f_{r}(x)B_{r}w\right) 
\quad\mathrm{mod}w,
$$
and so
\begin{equation}\label{v}
g_{r}(v_{(x,r)},v_{(x,r)})=4f_{r}(x)^{2}\left(\frac{(B_{r}^{2}w,w)}{(B_{r}w,w)}
-2r^{2}\gamma \right) ,
\end{equation}
where $w\in x$ is non-zero.
It follows that
\begin{equation}\label{svv}
S(V,V)=-2f\left(\frac{(B_{r}^{2}w,w)}{(B_{r}w,w)} -2r^{2}\gamma
\right) -4\lambda r^{2}f-6f^{2}.
\end{equation}
Relation (\ref{sp}), together with (\ref{svv}), determine $S(V).$
It remains to calculate $\mathrm{trace}_{\mathbb{R}}(S).$ Using
(\ref{s}) and (\ref{svv}), we have
\begin{align*}
\mathrm{trace}_{\mathbb{R}}(S)&=
\frac{1}{r^{2}}\left(\mathrm{trace}_{\mathbb{R}}(S_{r})-4mf\right)
+\frac{2}{r^{2}f}S(V,V)\\
&=-\frac{4}{r^{2}}\left(
(m+3)\left(\frac{(B_{r}^{2}w,w)}{B_{r}w,w)} +f\right)
-2r^{2}c\right) .
\end{align*}
Our claim follows now easily, combining the expressions of $S(L_{j})$, 
$S(v_{k})$, $S(V)$  and $\mathrm{trace}_{\mathbb{R}}(S)$ determined above.  

\end{proof}

We introduce a new family of hermitian operators
$\hat{B}_{r}:=\frac{1}{r^{2}}B_{r}=B+\delta (r) A$; $\delta
(r):=-r^{2}$ when $\lambda =0$ (equivalently, when $(M,\omega ,J)$
belongs to the third case of Proposition \ref{explicit}) and
$\delta (r)=-\frac{e^{\lambda r^{2}}}{\lambda }$ when $\lambda\neq
0$ (equivalently, when $(M,\omega ,J)$ belongs to the fourth case
of Proposition \ref{explicit}). Let $\hat{q}$, respectively
$\hat{Q}$, be the minimal and characteristic polynomials of
$\hat{B}_{r}$,  equal to
$$
\hat{q}(t)=(t-\gamma )^{3}\prod_{j=1}^{s}(t-\beta_{j});\quad
\hat{Q}:=(t-\gamma )^{n+1} \prod_{j=1}^{s}(t-\beta_{j})^{n_{j}}
$$
if there is $\mu_{j}\neq 0$ and
$$
\hat{q}(t)=(t-\gamma )^{2}\prod_{j=1}^{s}(t-\beta_{j});\quad
\hat{Q}(t)=(t-\gamma )^{n+1}\prod_{j=1}^{s}(t-\beta_{j})^{n_{j}}
$$
otherwise. Let $\hat{p}_{r}(t):=\prod_{j=1}^{l}(t-\eta_{j})$ be
the non-constant part of the characteristic polynomial
of the Bryant modified Ricci operator
$\hat{\Theta}_{r}$ of the K\"{a}hler structure $M_{\hat{B}_{r}}.$
It will be considered as a polynomial with function coefficients
defined on $\Sigma^{2m+1}_{\hat{B}_{r}}.$ We shall denote by 
$\hat{p}_{r,x}$ its value at a null line $x\in\Sigma^{2m+1}_{\hat{B}_{r}}$,
which is a polynomial with constant coefficients.

\begin{prop}
Let $\hat{q}_{1}$ be the constant polynomial equal to the quotient
of $\hat{q}$ by $(t-\gamma )^{2}$. Then the characteristic
polynomial $P(t)$ of the modified Ricci operator $\Theta$ of
$(M,g,J)$ is equal to the product $\frac{\hat{Q}\left(
t+\frac{c}{m+3}\right)} {\hat{q}\left( t+\frac{c}{m+3}\right)}
P_{1}(t)$, where
\begin{equation}\label{p1}
P_{1}(t):=
\left(t-\frac{ (m+2)c}{m+3}\right)\hat{p}_{r}
\left(t+\frac{c}{m+3}\right)
+\frac{f}{r^{2}}\hat{q}_{1}
\left( t+\frac{c}{m+3}\right).
\end{equation}
\end{prop}

\begin{proof} Using Proposition \ref{expresii},
together with the fact that the constant part of the
characteristic polynomial of the modified Ricci operator
$\Theta_{A}$ of a Bochner-flat K\"{a}hler structure $M_{A}$ is
equal to the quotient of the characteristic polynomial by the
minimal polynomial of $A$ (see \cite{paul}, Section 1.5), it is
easy to see that $P(t)=\frac{\hat{Q}\left( t+\frac{c}{m+3}\right)}
{\hat{q}\left( t+\frac{c}{m+3}\right)} P_{1}(t)$, where
\begin{align*}
P_{1}(t)&=\left( t-\frac{(m+2)c}{m+3}+\frac{f}{r^{2}}\right)
\hat{p}_{r}\left( t+\frac{c}{m+3}\right)\\
&+\frac{f}{r^{8}}
\sum_{j=1}^{l}\frac{q_{r}(\xi_{j})}{(\eta_{j}-\gamma
)^{2}p_{r}^{\prime}(\xi_{j})}\prod_{i\neq j}\left(
t-\eta_{i}+\frac{c}{m+3}\right) .
\end{align*}
We shall evaluate the expression
$$
\mathcal E (t):=\frac{f}{r^{8}}
\sum_{j=1}^{l}\frac{q_{r}(\xi_{j})}{(\eta_{j}-\gamma )^{2}p_{r}^{\prime}(\xi_{j})}\prod_{i\neq j}\left(
t-\eta_{i}+\frac{c}{m+3}\right) .
$$
For this, let $g_{r}, \hat{g}_{r}$ be the Bochner-flat K\"{a}hler
metrics of $M_{B_{r}}$, respectively, $M_{\hat{B}_{r}}$ (viewed as
metrics on $H$). Then, from relation (\ref{functg}), $\hat{g}_{r}=
r^{2}g_{r}=g$ on $H$. Using (\ref{grad}), we get
$$
-4\frac{q_{r}(\xi_{j})}{p_{r}^{\prime}(\xi_{j})}=
g_{r}\left(\mathrm{grad}_{g_{r}}(\xi_{j}),\mathrm{grad}_{g_{r}}(\xi_{j})\right)
=r^{6}\hat{g}_{r}
\left(\mathrm{grad}_{\hat{g}_{r}}(\eta_{j}),
\mathrm{grad}_{\hat{g}_{r}}(\eta_{j})\right)
=-4r^{6}\frac{\hat{q}(\eta_{j})}{\hat{p}_{r}^{\prime}(\eta_{j})},
$$
which implies that
$$
\mathcal E (t)=\frac{f}{r^{2}}\sum_{j=1}^{l}
\frac{\hat{q}_{1}(\eta_{j})}{\hat{p}_{r}^{\prime}(\eta_{j})}\prod_{i\neq
j} \left( t-\eta_{i}+\frac{c}{m+3}\right) .
$$
Note that $\mathcal E_{1}(t):=\frac{r^{2}}{f}{\mathcal E}(t)$
is a polynomial of degree $l-1$ which satisfies
$$
{\mathcal E}_{1}\left(\eta_{j}-\frac{c}{m+3}\right) =
\frac{\hat{q}_{1}(\eta_{j})}{\hat{p}_{r}^{\prime}(\eta_{j})}\prod_{i\neq j}
(\eta_{j}-\eta_{i})=\hat{q}_{1}(\eta_{j}),\quad j\in\{ 1,\cdots ,l\}
$$
Since $\hat{q}_{1}$ is a monic polynomial of degree $l$ and
${\mathcal E}_{1}$ is of degree $l-1$, we deduce that
$\hat{q}_{1}(t)={\mathcal E}_{1}\left( t-\frac{c}{m+3}\right)
+\hat{p}_{r}(t)$ which implies that
\begin{equation}\label{e}
\sum_{j=1}^{l}\frac{\hat{q}_{1}(\eta_{j})}{\hat{p}_{r}^{\prime}(\eta_{j})}\prod_{i\neq
j}\left( t-\eta_{i}+\frac{c}{m+3}\right) = \hat{q}_{1}\left(
t+\frac{c}{m+3}\right) - \hat{p}_{r}\left( t+\frac{c}{m+3}\right).
\end{equation}
Our claim follows.

\end{proof}

The possible constant roots of the polynomial $P_{1}$ (which are
also constant eigenvalues of $\Theta$) are determined in
Proposition \ref{constante}. In the proof of this Proposition we
will need the following additional Lemma.

\begin{lem}\label{dr}
The following equality holds:
$$
\frac{d}{dr}\hat{p}_{r}(t)=\frac{2f}{r}\left( \hat{p}_{r} \left(
t\right) - \hat{q}_{1}\left( t\right)\right).
$$
\end{lem}

\begin{proof}
We take the derivative with respect to $r$ of the equality
$$
\left( (\eta_{j}I-B-\delta A)^{-1}w,w\right) =0
$$
(which follows from $\hat{p}_{r}(\eta_{j})=0$) and we obtain
\begin{equation}\label{eta}
\left( ( \eta_{j}I-\hat{B}_{r})^{-1}(\dot{\eta}_{j}I-\dot{\delta}
A) (\eta_{j}I-\hat{B}_{r})^{-1}w,w\right)=0.
\end{equation}
On the other hand, since $AB=\gamma A$, it is easy to see that
\begin{equation}\label{eta1}
(\eta_{j}I-\hat{B}_{r})^{-1}A(\eta_{j}I-\hat{B}_{r})^{-1}=
\frac{A}{(\eta_{j}-\gamma )^{2}}.
\end{equation}
Applying (\ref{patrat}) to $A:=\hat{B}_{r}$ and using the
fact that $\eta_{j}I-\hat{B}_{r}$ is invertible, we get
\begin{equation}\label{dd}
\left( (\eta_{j}I-\hat{B}_{r})^{-2}w,w\right) =
-\frac{\hat{p}_{r}^{\prime}(\eta_{j})}{\hat{q}(\eta_{j})}
(\hat{B}_{r}w,w) .
\end{equation}
Combining (\ref{eta}), (\ref{eta1}) and (\ref{dd}), and using the
fact that $\delta (r)=-\frac{e^{\lambda r^{2}}}{\lambda}$ when
$\lambda\neq 0$ and $\delta (r)=-r^{2}$ when $\lambda =0$, we
deduce the expressions of the derivatives $\dot{\eta}_{j}$ as
follows:
\begin{equation}\label{dot}
\dot{\eta}_{j}=
\frac{2f}{r}\frac{\hat{q}_{1}(\eta_{j})}{\hat{p}_{r}^{\prime}(\eta_{j})}.
\end{equation}
Since
$\frac{d}{dr}\hat{p}_{r}(t)=-\sum_{j=1}^{l}\dot{\eta}_{j}\prod_{i\neq
j}(t-\eta_{i})$ we get, using (\ref{e}),  our claim.
\end{proof}

\begin{prop}\label{constante}  The following statements hold:

\begin{enumerate}

\item Suppose that $\hat{q}_{1}(c)\neq 0.$ Then
the polynomial $P_{1}$ defined in (\ref{p1}) does not have
constant roots, except when $\lambda\neq 0$ and $\alpha
=\mu_{j}=0$ for any $j$. In this case,
$t:=\frac{(m+2)c}{m+3}+\lambda$ is the unique constant root of $P_{1}$ and is
simple.

\item Suppose that $\hat{q}_{1}(c)=0.$ Then $t:=\frac{(m+2)c}{m+3}$ is
a simple root of $P_{1}$. The polynomial $P_{1}$ has other
constant roots if and only if $\lambda\neq 0$ and $\alpha
=\mu_{j}=0$ for any $j$. In this case, there is only one
additional constant root of $P_{1}$, namely $t:=\lambda +\frac{(m+2)c}{m+3}$,
which is simple.

\end{enumerate}
\end{prop}

\begin{proof} We
first consider the case when $\hat{q}_{1}(c)\neq 0.$ We claim that
$P_{1}$ has no multiple roots. Suppose, on the contrary, that $t$
is a multiple (necessarily constant, because the non-constant
eigenvalues of the Bryant modified Ricci operator are always simple) root
of $P_{1}.$ Since $\hat{q}_{1}(c)\neq 0$,  $t$ cannot be equal to
$\frac{(m+2)c}{m+3}$ and so $\hat{q}_{1}\left(
t+\frac{c}{m+3}\right)\neq 0$ (because 
$\hat{p}_{r}$ has no
constant roots). The equalities $P_{1}(t)=P_{1}^{\prime}(t)=0$
imply that
\begin{equation}\label{cannot}
t_{1}\hat{p}_{r}^{\prime}\left(t+\frac{c}{m+3}\right) +
t_{2}\hat{p}_{r}\left(t+\frac{c}{m+3}\right)=0,
\end{equation}
where $t_{1}:= t-\frac{(m+2)c}{m+3}\in\mathbb{R}\setminus \{ 0\}$
and $t_{2}:= -\frac{(m+2)c}{m+3}-\left(
t-\frac{(m+2)c}{m+3}\right)
\frac{\hat{q}_{1}^{\prime}\left(t+\frac{c}{m+3}\right)}
{\hat{q}_{1}^{\prime}\left(t+\frac{c}{m+3}\right)}\in\mathbb{R} .$
But (\ref{cannot}) cannot hold: if it did, it would imply
that $I,\hat{B}_{r},\cdots ,\hat{B}_{r}^{l+1}$ were dependent,
which contradicts the fact that the minimal polynomial of
$\hat{B}_{r}$ has degree $l+2$. We conclude that $P_{1}$ cannot
have multiple roots. We will now show that the only possible
constant root of $P_{1}$ is $\lambda +\frac{(m+2)c}{m+3}$ and it
is a root if and only if $\lambda\neq 0$ and $\alpha =\mu_{j}=0$
for any $j$. For this, let $t$ be a constant root of $P_{1}.$
Taking the derivative with respect to $r$ of the equality
$P_{1}(t)=0$ and using Lemma \ref{dr}, we get
$$
\left( t-\frac{(m+2)c}{m+3}-\lambda\right)\hat{q}_{1}\left(
t+\frac{c}{m+3}\right) =0,
$$
from where we deduce that $t=\lambda +\frac{(m+2)c}{m+3}$, since
$\hat{q}_{1}\left(t+\frac{c}{m+3}\right)\neq 0$. Moreover,
$P_{1}(t)=0$ if and only if
\begin{equation}\label{ga}
\lambda\hat{p}_{r}(\gamma )+\frac{f}{r^{2}}\hat{q}_{1}(\gamma )=0.
\end{equation}
Equality (\ref{ga}) forces $\lambda\neq 0$ (if $\lambda =0$, then,
from (\ref{ga}), $\hat{q}_{1}(\gamma )=0$; also, $c=\gamma$;
recall however that we are under the hypothesis
$\hat{q}_{1}(c)\neq 0$; we obtain a contradiction). Therefore,
$\lambda\neq 0$ and then $A_{r}=r^{4}e^{\lambda r^{2}}A$. Relation
(\ref{ga}) is equivalent with $\lambda \left(
\widetilde{\hat{b}_{r}}(\gamma )w,w\right) + \hat{q}_{1}(\gamma
)e^{\lambda r^{2}}(Aw,w)=0$,  for any $w\in W$ null, where
$\widetilde{\hat{b}}_{r}$ denotes the reduced adjoint operator of
$\hat{B}_{r}.$ With the notations of Remark \ref{op},
$\widetilde{\hat{b}_{r}}(\gamma )$, as well as $A$, act trivially
on the subspaces $W_{j}$ (for $j\geq 1$) of $W$ (here and below
the Reader is referred to Lemma $3$ of \cite{paul}, which
describes the action of the reduced adjoint operator of a $k$-step
parabolic hermitian operator when applied to the parabolic
eigenvalue). It follows that (\ref{ga}) is equivalent with
\begin{equation}\label{ga1}
\lambda\widetilde{\hat{b}_{r}}(\gamma ) + \hat{q}_{1}(\gamma
)e^{\lambda r^{2}}A=0.
\end{equation}
We claim that equality (\ref{ga1}) holds if and only if $\alpha
=\mu_{j}=0$ for any $j$ (and $\lambda\neq 0$). Notice first that
if (\ref{ga1}) holds then $\hat{q}_{1}(\gamma )\neq 0$ (since
$\widetilde{\hat{b}_{r}}(\gamma )\neq 0$), which implies that
$\hat{B}_{r}$ is $1$-step parabolic. We deduce that $\mu_{j}=0$
for any $j.$ With the notations of Remark \ref{op},
$\widetilde{\hat{b}}_{r}(\gamma )$, as well as $A$, act trivially
on $\widehat{{W}_{0}}$; on $\mathrm{Span}\{ e_{0},e_{1}\}$,
$\widetilde{\hat{b}_{r}}(\gamma )$ acts by $\lambda\left( \alpha
-\frac{e^{\lambda r^{2}}}{\lambda}\right)\prod_{j=1}^{s}(\gamma
-\beta_{j})A$. Relation (\ref{ga1}) is equivalent to $\left(
\alpha\lambda -e^{\lambda r^{2}}\right)\prod_{j=1}^{s} (\gamma
-\beta_{j})+\hat{q}_{1}(\gamma )e^{\lambda r^{2}}=0$; since
$\hat{q}_{1}(\gamma )=\prod_{j=1}^{s} (\gamma -\beta_{j})$ and $\lambda\neq 0$, it
reduces to $\alpha =0.$
Since $P_{1}$ doesn't have multiple roots, $\lambda+\frac{(m+2)c}{m+3}$
is the (unique) simple root of $P_{1}.$\\

Now we consider the case when $\hat{q}_{1}(c)=0.$ Clearly,
$\frac{(m+2)c}{m+3}$ is a root in this case. Define the polynomial
$\hat{q}_{2}(t):=\frac{\hat{q}_{1}(t)}{t-c}$ and
let $t$ be a
constant root of 
\begin{equation}\label{newec}
\hat{p}_{r}\left( t+\frac{c}{m+3}\right)
+\frac{f}{r^{2}}\hat{q}_{2} \left( t+\frac{c}{m+3}\right) =0.
\end{equation}
Taking the derivative with respect to $r$ of (\ref{newec})
and using Lemma \ref{dr}, we get
$$
\left( t-\frac{(m+2)c}{m+3}-\lambda \right)\hat{q}_{2}\left(
t+\frac{c}{m+3}\right) = 0.
$$
This implies that $t =\lambda +\frac{(m+2)c}{m+3}$, because
$\hat{q}_{2}\left( t+\frac{c}{m+3}\right)\neq 0$ (because
$\hat{p}_{r}$ has no constant roots). As before,  $t$
is a root of (\ref{newec}) if and only if
\begin{equation}\label{ec1}
\widetilde{\hat{b}_{r}}(\gamma ) + \hat{q}_{2}(\gamma )e^{\lambda
r^{2}}A=0.
\end{equation}
Notice first that if (\ref{ec1}) holds, then $\hat{q}_{2}(\gamma
)\neq 0$. Next, we prove that if (\ref{ec1}) holds then
$\hat{B}_{r}$ is $1$-step parabolic.
The argument is the following: 
suppose, on the contrary,
that (\ref{ec1}) holds and that $\hat{B}_{r}$ is $2$-step
parabolic; then $\lambda = 0$; otherwise, since $(t-\gamma )$
divides $\hat{q}_{1}$ ($\hat{B}_{r}$ being $2$-step parabolic),
$\hat{q}_{2}(\gamma )=0$, which is impossible. 
On the other hand,
$\widetilde{\hat{b}_{r}}(\gamma )$ acts as
$\sum_{k=2}^{n}|\mu_{k}|^{2}\prod_{j=1}^{s}(\gamma -\beta_{j})A$
on $\mathrm{Span}\{ e_{0},e_{1}\}$,
when $\hat{B}_{r}$ is $2$-step parabolic; also $\hat{q}_{2}(\gamma
)=\prod_{j=1}^{s}(\gamma -\beta_{j})$, when $\lambda =0;$ from (\ref{ec1}) it follows
that 
$$
\left(\sum_{k=2}^{n}|\mu_{k}|^{2}+1\right)\prod_{j=1}^{s}(\gamma -\beta_{j})A=0
$$ 
on $\mathrm{Span}\{ e_{0},e_{1}\}$, 
which cannot hold. We have proved that 
if $\lambda +\frac{(m+2)c}{m+3}$ is a root, then
$\hat{B}_{r}$ is
$1$-step parabolic. Moreover, in this case $\lambda\neq 0$ (if $\lambda =0$ then $c
=\gamma$ and since $\hat{q}_{1}(c )=0$, then $\hat{q}_{1}(\gamma
)=0$ which is absurd because $(t-\gamma )^{3}$ does not divide the
minimal polynomial $\hat{q}$ of $\hat{B}_{r}$ when $\hat{B}_{r}$
is $1$-step parabolic). Finally, when $\hat{B}_{r}$ is $1$-step
parabolic and $\lambda \neq 0$, relation (\ref{ec1}) becomes
$$
\left(\alpha -\frac{e^{\lambda r^{2}}}{\lambda}\right)
\prod_{j=1}^{s}\left(\gamma -\beta_{j}\right) +\hat{q}_{2}(\gamma
)e^{\lambda r^{2}}=0
$$
which holds if and only if $\alpha =0$, because
$$
\hat{q}_{2}(\gamma )=\frac{1}{\gamma -c}\hat{q}_{1}(\gamma )=
\frac{1}{\lambda }\hat{q}_{1}(\gamma )=
\frac{1}{\lambda}\prod_{j=1}^{s}(\gamma -\beta_{j}).
$$
We have proved that $P_{1}$ has an additional constant root, besides
$\frac{(m+2)c}{m+3}$, if and only if $\hat{B}_{r}$ is $1$-step parabolic 
(i.e. $\mu_{j}=0$ for any $j$), $\lambda\neq 0$ and $\alpha =0.$ The additional
constant root is $\lambda +\frac{(m+2)c}{m+3}$.   
It is easy to see that it is simple.

\end{proof}

In Proposition \ref{2st} we shall determine the Bryant minimal and
characteristic polynomials of $(M,\omega ,J)$. For the proof of
this Proposition we need the following additional Lemma.

\begin{lem}\label{hori} For any $t\in\mathbb{R}$,
$(x,r)\in M$ and $w\in x$ non-zero,
\begin{align*}
g_{(x,r)}\left( d^{H}\hat{p}_{r}(t),d^{H}\hat{p}_{r}(t)\right) &= 4\left(
\hat{q}^{\prime}(t)\hat{p}_{r,x}(t)-\hat{q}(t)\hat{p}_{r,x}^{\prime}(t)
-2t\hat{p}_{r,x}^{2}(t)
+\hat{p}_{r,x}^{2}(t)\frac{(B_{r}^{2}w,w)}{r^{2}(B_{r}w,w)}
\right)\\
g_{(x,r)}\left(d^{H}\left(
\frac{f}{r^{2}}\right),d^{H}\left(\frac{f}{r^{2}}\right) \right)
&=\frac{4f^{2}}{r^{6}}\left(\frac{(B_{r}^{2}w,w)}{(B_{r}w,w)}
-2r^{2}\gamma \right)\\
g_{(x,r)}\left(d^{H}\left(
\frac{f}{r^{2}}\right),d^{H}\hat{p}_{r}(t)\right) &=
\frac{4f}{r^{2}}\left( (t-\gamma )\hat{q}_{1}(t)-(t+
\gamma)\hat{p}_{r,x}(t)
+\frac{\hat{p}_{r,x}(t)(B_{r}^{2}w,w)}{r^{2}(B_{r}w,w)}\right) .\\
\end{align*}
\end{lem}

\begin{proof}
The first equality follows from Lemma \ref{aditional}, applied to
$\hat{B}_{r}$. To prove the second equality, we remark that the
$1$-form $d^{H}\left( \frac{f}{r^{2}}\right)$ corresponds to the
vector field $\frac{1}{r^{4}}v$ by means of the metric $g$.
Therefore, $g\left(d^{H}\left(
\frac{f}{r^{2}}\right),d^{H}\left(\frac{f}{r^{2}}\right)\right)$
is equal to $\frac{1}{r^{6}}g_{r}(v,v).$ The second equality of
the Lemma follows from (\ref{v}). To prove the third
equality we notice that $g\left(d^{H}\left(
\frac{f}{r^{2}}\right),d^{H}\hat{p}_{r}(t)\right)$ is equal to
$\frac{\hat{L}_{t}(f)}{r^{2}},$ where $\hat{L}_{t}$ is the vector
field on $M$ which, at a point $(x,r)\in M$, belongs to 
$H_{x}=\mathrm{Hom}_{\mathbb{C}}(x,W/x)$ and is the endomorphism
$$
(\hat{L}_{t})_{(x,r)}(w)= 2\left(
\tilde{\hat{b}}_{r}(t)w-\hat{p}_{r}(t)\hat{B}_{r}w\right) \quad
\mathrm{mod}w,\quad w\in x.
$$
This is true since $\hat{L}_{t}$ corresponds to the $1$-form 
$d^{H}\hat{p}_{r}(t)$ by means
of the metric $g$ (which coincides with $\hat{g}_{r}$ on the
bundle $H$, restricted to a level set $N_{r}$). It is
straightforward to check  that
\begin{align*}
(A_{r}\hat{L}_{t}w,w)&=2\left( (t-\gamma )
\hat{q}_{1}(t)-\gamma \hat{p}_{r,x}(t)\right) (A_{r}w,w)\\
(B_{r}\hat{L}_{t}w,w)&=2\hat{p}_{r,x}(t)\left( t-
\frac{(B_{r}^{2}w,w)}{r^{2}(B_{r}w,w)}\right) (B_{r}w,w),
\end{align*}
so that
$$
\hat{L}_{t}(f)_{(x,r)}=4f\left( \left(t-\gamma \right)\hat{q}_{1}(t)-(t+
\gamma) \hat{p}_{r,x}(t)
+\frac{\hat{p}_{r,x}(t)(B_{r}^{2}w,w)}{r^{2}(B_{r}w,w)}\right) ,
$$
which proves the third equality.

\end{proof}

\begin{prop}\label{2st}
The Bryant characteristic polynomial of $(M,g,J)$ is equal to
$$
p_{c}(t)=\left( t-\frac{(m+2)c}{m+3}\right) \left(
t+\frac{c}{m+3}-\gamma\right)^{n+1}\prod_{j=1}^{s} \left(
t+\frac{c}{m+3}-\beta_{i}\right)^{n_{j}}
$$
The Bryant minimal polynomial $p_{m}$ of $(M,g,J)$ has the
following expression:
\begin{enumerate}
\item if there is $\mu_{j}\neq 0$
and $c$ is different from $\beta_{k}$ (for any $k$) and $\gamma$,
then
$$
p_{m}(t)=\left( t-\frac{(m+2)c}{m+3}\right) \left(
t+\frac{c}{m+3}-\gamma\right)^{3}\prod_{j=1}^{s}\left(
t+\frac{c}{m+3}-\beta_{j}\right) .
$$

\item if there is $\mu_{j}\neq 0$ and $c$ is equal to $\beta_{j}$
(for a certain $j$) or to $\gamma$, then
$$
p_{m}(t)= \left( t+\frac{c}{m+3}-\gamma\right)^{3}
\prod_{j=1}^{s}\left( t+\frac{c}{m+3}-\beta_{j}\right) .
$$

\item if all $\mu_{j}=0$ and $c$ is different from $\beta_{j}$ (for any $j$)
then
$$
p_{m}(t)=\left( t-\frac{(m+2)c}{m+3}\right) \left(
t+\frac{c}{m+3}-\gamma\right)^{2}\prod_{j=1}^{s}\left(
t+\frac{c}{m+3}-\beta_{j}\right).
$$
except when $\lambda\neq 0$ and $\alpha =0$, in which case
$$
p_{m}(t)=\left( t-\frac{(m+2)c}{m+3}\right) \left(
t+\frac{c}{m+3}-\gamma\right)\prod_{j=1}^{s}\left(
t+\frac{c}{m+3}-\beta_{j}\right).
$$

\item if all $\mu_{j}=0$ and $c$ is equal to $\beta_{j}$ (for a certain $j$),
then
$$
p_{m}(t)= \left(
t+\frac{c}{m+3}-\gamma\right)^{2}\prod_{j=1}^{s}\left(
t+\frac{c}{m+3}-\beta_{j}\right)
$$
except when $\alpha =0$ (and $\lambda\neq 0$), when
$$
p_{m}(t)= \left(
t+\frac{c}{m+3}-\gamma\right)\prod_{j=1}^{s}\left(
t+\frac{c}{m+3}-\beta_{j}\right) .
$$
\end{enumerate}
\end{prop}

\begin{proof}
Let $t$ be a non-constant root of the polynomial $P_{1}$. Then
\begin{align*}
\left( t-\frac{(m+2)c}{m+3}\right)\hat{p}_{r}\left(
t+\frac{c}{m+3}\right)
&=-\frac{f}{r^{2}}\hat{q}_{1}\left( t+\frac{c}{m+3}\right)\\
\left( t-\frac{(m+2)c}{m+3}\right)^{2}
\hat{p}_{r}^{\prime}\left(t+\frac{c}{m+3}\right)&= \left(
t-\frac{(m+2)c}{m+3}\right) P_{1}^{\prime}(t)+\frac{f}{r^{2}}
\hat{q}_{1}\left( t+\frac{c}{m+3}\right)\\
&-\frac{f}{r^{2}}\left( t-\frac{(m+2)c}{m+3}\right)
\hat{q}_{1}^{\prime} \left( t+\frac{c}{m+3}\right) .
\end{align*}
Using these relations and Lemma \ref{hori}, we can calculate the
square norm of the covector $(d^{H}P_{1})(t)$ as follows:
\begin{align*}
g\left( (d^{H}P_{1})(t),(d^{H}P_{1})(t)\right)&= -4\hat{q}\left(
t+\frac{c}{m+3}\right)\left( t-\frac{(m+2)c}{m+3}\right)
P_{1}^{\prime}(t)\\
&-\frac{4f}{r^{2}}\left( t-\frac{(m+2)c}{m+3}\right) \hat{q}^{\prime}\left(
t+\frac{c}{m+3}\right)\hat{q}_{1}\left( t+\frac{c}{m+3}\right)\\
&-\frac{4f}{r^{2}}\hat{q}\left( t+\frac{c}{m+3}\right)\hat{q}_{1}
\left( t+\frac{c}{m+3}\right)\\
&+\frac{4f}{r^{2}}
\left( t-\frac{(m+2)c}{m+3}\right) \hat{q}\left( t+\frac{c}{m+3}\right)
\hat{q}_{1}^{\prime}\left( t+\frac{c}{m+3}\right)\\
&+\frac{8f}{r^{2}}\left( t-\frac{(m+2)c}{m+3}\right)\left( t+
\frac{c}{m+3}-\gamma \right)\hat{q}_{1}^{2}\left( t+\frac{c}{m+3}\right).
\end{align*}
Since
$\hat{q}^{\prime}(t)\hat{q}_{1}(t)-\hat{q}(t)\hat{q}_{1}^{\prime}(t)=
2(t-\gamma )\hat{q}_{1}^{2}(t),$ the above expression reduces to
\begin{align*}
g((d^{H}P_{1})(t),(d^{H}P_{1})(t))&=-4\hat{q}\left(
t+\frac{c}{m+3}\right)
\left( t-\frac{(m+2)c}{m+3}\right) P_{1}^{\prime}(t)\\
&-\frac{4f}{r^{2}}\hat{q}_{1}^{2}\left(
t+\frac{c}{m+3}\right)\left( t-\gamma +\frac{c}{m+3}\right)^{2}.
\end{align*}

On the other hand, Lemma \ref{dr} together with the definition of
$P_{1}$ imply that
$$
\left(\frac{d}{dr}P_{1}\right) (t)=-\frac{2f}{r}\hat{q}_{1}\left(
t+\frac{c}{m+3}\right) \left(
t-\frac{(m+2)c}{(m+3)}-\lambda\right) .
$$
Since $g(dr, dr)=\frac{1}{f}$, we obtain
\begin{align*}
g\left( (dP_{1})(t),(dP_{1})(t)\right) &=g\left(
(d^{H}P_{1})(t),(d^{H}P_{1})(t)\right)
+\left(\frac{d}{dr}P_{1}(t)\right)^{2}g\left( dr,dr\right)\\
&=-4\hat{q}\left( t+\frac{c}{m+3}\right)
\left( t-\frac{(m+2)c}{m+3}\right) P^{\prime}_{1}(t).
\end{align*}
On the other hand, since $P_{1}^{\prime}(t)dt+d\left(
P_{1}(t)\right) =0$, we get
\begin{equation}\label{dt}
g(dt,dt )=\frac{g\left( (dP_{1})(t) ,(d P_{1})(t)\right)}{
P_{1}^{\prime}(t)^{2}}.
\end{equation}
We distinguish three cases: (i) $P_{1}$ has no constant roots; (ii)
$P_{1}$ has a unique constant root, which is simple and
equal to 
$t_{1}=\lambda +\frac{(m+2)c}{m+3}$; (iii) 
$P_{1}$ has two constant roots, $t_{1}$ and $t_{2}=
\frac{(m+2)c}{m+3}$, which are simple and
distinct (see
Proposition \ref{constante}). In all cases, the Bryant
characteristic polynomial of $(M,\omega ,J)$ is
$$
p_{c}(t)=\left( t-\frac{(m+2)c}{m+3}\right) \hat{Q}\left(
t+\frac{c}{m+3}\right).
$$
In case (i), the Bryant minimal polynomial of $(M,g,J)$ is
$$
p_{m}(t)=\hat{q}\left( t+\frac{c}{m+3}\right) \left(
t-\frac{(m+2)c}{m+3}\right)
$$
In case (ii), 
$$
p_{m}(t)= \left( t-t_{1}\right)^{-1}\left(
t-\frac{(m+2)c}{m+3}\right) \hat{q}\left(t +\frac{c}{m+3}\right).
$$
and in case (iii),
$$
p_{m}(t)= \left( t-t_{1}\right)^{-1}
\left( t-t_{2}\right)^{-1}
\left( t-\frac{(m+2)c}{m+3}\right) \hat{q}\left(t +\frac{c}{m+3}\right).
$$

\end{proof}

\textbf{Proof of Proposition \ref{es}:} Proposition \ref{es}
is an easy consquence Proposition \ref{2st}.

\section{Examples}

In this section we consider some important classes of Bochner-flat
K\"{a}hler manifolds and we show how they can be realised locally
as generalised K\"{a}hler cones.\\

(i) {\it Bryant's Bochner-flat K\"{a}hler structures.} Let
$N=S^{2m+1}\subset\mathbb{C}^{m+1}$ with its standard CR structure
and $(k_{1},\cdots ,k_{m+1})$ a system of non-negative real
numbers. Define, for every $r>0$, the vector field
$$
T_{r}(z):=\sum_{j=1}^{m+1}\left(1+k_{j}r^{2}\right)\left(
x_{j}\frac{\partial}{\partial y_{j}}-
y_{j}\frac{\partial}{\partial x_{j}}\right) ,
$$
which is the Reeb vector field of a Sasaki structure on
$S^{2m+1}\subset\mathbb{C}^{m+1}.$ Here $z=(z_{1},\cdots ,z_{m+1})$ belongs to
$S^{2m+1}$, $z_{j}=x_{j}+iy_{j}$ for any $j\in \{ 1,\cdots ,m+1\}
$ and $r^{2}=|z_{1}|^{2}+\cdots +|z_{m+1}|^{2}.$ The family of
Sasaki Reeb vector fields $\{ T_{r},r>0\}$ defines a Bochner-flat
generalised K\"{a}hler cone on $\mathbb{C}^{m+1}\setminus \{ 0\}$,
which belongs to the first case of Proposition \ref{explicit}; the
solution of equation (\ref{muec}) is $\mu (t)=-\frac{2}{t}$ and
the hermitian operator $B$ is semi-simple, with eigenvalues
$k_{j}^{\prime}=k_{j}-\frac{1}{m+2}\sum_{i=1}^{m+1}k_{i}$, for
$j\in \{ 1,\cdots ,m+1\} .$ This Bochner-flat K\"{a}hler structure
has been discovered by Robert Bryant in \cite{bryant} and has been
further studied in \cite{paul}; it can be extended as a complete
Bochner-flat K\"{a}hler structure on
$\mathbb{C}^{m+1}.$\\

(ii) {\it Bochner-flat K\"{a}hler-Einstein structures}. Let
$(M,g,J)$ be a Bochner-flat generalised K\"{a}hler cone. With the
notations of Proposition \ref{explicit}, suppose that $B=eA$, for
$e\in\mathbb{R}.$ If $(M,g,J)$ belongs to the first and second
case of Proposition \ref{explicit}, then it is K\"{a}hler-Einstein
if and only if $e^{2}+2d=0$; moreover, the Bryant modified Ricci
operator of $(M,\omega ,J)$ is $\Theta =\frac{e}{m+3}\mathrm{Id}$
in the first case and $\Theta =-\frac{e}{m+3}\mathrm{Id}$ in the
second case. If $(M,g,J)$ belongs to the third case of Proposition
\ref{explicit} and $B=eA$ then it is never Einstein; if it belongs
to the fourth case of Proposition \ref{explicit} then it is
K\"{a}hler-Einstein if and only if $e=0$ (and $\lambda <0$); The
Bryant modified Ricci operator is $\Theta
=\frac{\lambda}{m+3}\mathrm{Id}$ in this
case.\\

(iii) {\it Bochner-flat generalised K\"{a}hler cones of order
one.} If $B=eA$ but $(M,g,J)$ is not Einstein, then it must have
order one. The Bryant minimal and characteristic polynomials have
the following expressions: if $(M,g,J)$ belongs to the case $1$,
respectively to the case $2$ of Proposition \ref{explicit}, then
\begin{align*}
p_{c}(t)&=\left( t\mp \frac{e}{m+3}\right)^{m+1} \left( \left(
t\mp\frac{e}{m+3}\right)^{2}\pm
e\left( t\mp\frac{e}{m+3}\right) +\frac{e^{2}+2d}{4}\right)\\
p_{m}(t)&=\left( t\mp \frac{e}{m+3}\right) \left( \left(
t\mp\frac{e}{m+3}\right)^{2}\pm
e\left( t\mp\frac{e}{m+3}\right)+\frac{e^{2}+2d}{4}\right) ;\\
\end{align*}
$(M,g,J)$ is of hyperbolic type when $d>0$, of $1$-step parabolic type
when $d=0$ and of elliptic type when $d<0.$ If $(M,J,g )$ belongs
to the cases $3$ and $4$ of Proposition \ref{explicit}, then
\begin{align*}
p_{c}(t)&=\left( t+\frac{(m+2)\lambda}{m+3}\right)\left(
t-\frac{\lambda}{m+3}\right)^{m+2}\\
p_{m}(t)&=\left( t+\frac{(m+2)\lambda}{m+3}\right)\left(
t-\frac{\lambda}{m+3}\right)^{2}.
\end{align*}
(Recall that $\lambda =0$ in the case $3$ and $\lambda\neq 0$ in
case $4$). When $(M,\omega ,J)$ belongs to the case $3$, it is
$2$-step parabolic; when it belongs to case $4$, it is $1$-step
parabolic. Bochner-flat K\"{a}hler structures of order one have
been studied in \cite{acg}. As shown in \cite{acg}, \cite{paul},
they fiber over a K\"{a}hler manifold with constant holomorphic
sectional curvature (in our formalism, the fibration is
$M\to
{N}_{r}/\tilde{T}_{r}$).\\

(iv) {\it Weighted projective spaces as generalised K\"{a}hler
cones.} Let $\mathbb{P}_{a}^{m+1}$ be a weighted projective space,
of weights $(a_{1},\cdots ,a_{m+1})$, where $a_{j}$ are positive
integers. As shown in \cite{bryant}, \cite{paul},
$\mathbb{P}_{a}^{m+1}$ has a canonical Bochner-flat K\"{a}hler
structure, of semi-simple type, isomorphic with $M_{C}$, where $C$
is a hermitian semi-simple operator of $\mathbb{C}^{m+2,1}$, with
eigenvalues $-\sum_{j=1}^{m+2}\lambda_{j},\lambda_{1},\cdots
,\lambda_{m+2}$, where $\lambda_{j}$ are related to the weights
$a_{j}$ by $\lambda_{j}=a_{j}-\frac{1}{m+3}\sum_{i=1}^{m+2}a_{i},$
for any $j\in \{ 1,\cdots ,m+2\}. $ As a Bochner-flat generalised
K\"{a}hler cone, $\mathbb{P}_{a}^{m+1}$ belongs to the first case
of Proposition \ref{explicit}; $\mu$ is any solution of equation
(\ref{muec}), with $d=\frac{2a_{m+2}^{2}}{(m+3)^{2}}$, and the
hermitian operator $B$ is semi-simple, with eigenvalues
$-\frac{\sum_{j=1}^{m+1}a_{j}}{m+2},a_{1}-\frac{\sum_{j=1}^{m+1}a_{j}}{m+2},
\cdots ,a_{m+1}-\frac{\sum_{j=1}^{m+1}a_{j}}{m+2}$.\\

(v) {\it Tachibana and Liu Bochner-flat K\"{a}hler generalised
K\"{a}hler cones.} Consider a Bochner-flat generalised K\"{a}hler
cone structure $(g,J)$ which belongs to the first case of
Proposition \ref{explicit}. With the notations of Section
\ref{quotient}, assume that $k_{1}=\cdots =k_{m+1}:=\bar{k}$. On
the set $D\subset\mathbb{C}^{m+1}$ defined by (\ref{d}), the
K\"{a}hler structure $(g,J)$ has a K\"{a}hler potential $x$ which
depends only on $r^{2}:=|z_{1}|^{2}+ \cdots +|z_{m+1}|^{2}$. As a
function of $r^{2}=t$, $x$ satisfies the implicit equation
\begin{equation}\label{ec}
e^{ax(t)}\dot{\mu}(x(t))^{-\frac{1}{2}}=t,
\end{equation}
where $a:=(m+2)\bar{k}.$ In general, a K\"{a}hler structure which
is defined on an open subset of the standard $\mathbb{C}^{m}$ and
has a global K\"{a}hler potential, say $h$, which depends only on
$r^{2}$ is Bochner-flat, if and only if $h$, as a function of
$r^{2}=t$, satisfies a differential equation of the form  \cite{liu}
\begin{equation}\label{tls}
\ddot{h}(t)=\lambda_{1}t\dot{h}^{3}(t)+\lambda_{2}\dot{h}^{2}(t),
\end{equation}
where $\lambda_{1},\lambda_{2}\in\mathbb{R}$. It can be easily
verified that if $x$ satisfies (\ref{ec}), then it satisfies also
(\ref{tls}), with $\lambda_{1}:=a^{2}+\frac{d}{2}$ and
$\lambda_{2}:=-2a.$

\vspace{10pt}

\textbf{Author's address:}
{\it Present:} Centro di Ricerca Matematica Ennio de Giorgi, Scuola Normale
Superiore di Pisa, Collegio Puteano, Piazza dei Cavalieri nr. 3, Pisa,
Italia.

{\it Permanent:}  Institute of Mathematics of the
Romanian Academy "Simion Stoilow", Calea Grivitei nr. 21, Sector 1,
Bucharest, Romania; e-mail: liana.david@imar.ro ; lili@mail.dnttm.ro

\end{document}